\documentclass[11pt]{article}
\usepackage[utf8]{inputenc}

\usepackage{geometry}
\usepackage{amsmath, amsfonts}
\usepackage[noend]{algorithmic}
\usepackage{algorithm}
\usepackage{color}
\usepackage{natbib}
\usepackage{lscape}
\usepackage{graphicx}
\usepackage{caption}
\usepackage{subcaption}
\usepackage{multirow}
\usepackage{booktabs}
\usepackage{bm}

\usepackage{rotating}

\usepackage[colorinlistoftodos]{todonotes}

\newcommand{\proof}{\noindent{\bf Proof}\hspace{0.2cm}}
\newcommand{\prend}{\hfill\fbox{}\vskip 6pt}

\usepackage{xurl}

\newcommand{\spmodel}{{\ttfamily{SP-CFL}}}
\newcommand{\mpmodel}{{\ttfamily{MP-CFL}}}

\newtheorem{theorem}{Theorem}

\title{\bf Incorporating time-dependent demand patterns in the optimal location of capacitated charging stations \\}

\author{ \it Carlo Filippi$ ^{(1)}$ \quad Gianfranco Guastaroba$ ^{(1)}$\\ \it Lorenzo Peirano$ ^{(1)}$ \quad M. Grazia Speranza$ ^{(1)}$\\
{\small $^{(1)}$ \it University of Brescia, Department of Economics and Management, Brescia, Italy}\\
{\small \{carlo.filippi, gianfranco.guastaroba, lorenzo.peirano, grazia.speranza\}@unibs.it}\\
}
\date{\today}

\begin{document}

\maketitle \thispagestyle{empty}

\begin{abstract}

A massive use of electric vehicles is nowadays considered to be a key element of a sustainable transportation policy and the availability of charging 
stations is a crucial issue for their extensive use. Charging stations in an urban area have to be deployed in such a way that they can satisfy a demand 
that may dramatically vary in space and time. In this paper we present an optimization model for the location of charging stations that takes into 
account the main specific features of the problem, in particular the different charging technologies, and their associated service time, and the fact 
that the demand depends on space and time. To measure the importance of incorporating the time dependence in an optimization model, we also present a simpler 
model that extends a classical location model and does not include the temporal dimension. A worst-case analysis and extensive computational experiments 
show that ignoring the temporal dimension of the problem may lead to a substantial amount of unsatisfied demand.

\end{abstract}

\noindent
{\emph{Keywords:} Facility location, Charging stations, Electric vehicles, Demand patterns, Time-dependent optimization.}

\section{Introduction}

Sustainable transportation is one of the major challenges that modern countries are facing. Several sources indicate that the transportation sector 
generates the largest share of GreenHouse Gas (GHG) emissions. According to the United States Environmental Protection 
Agency\footnote{\url{https://www.epa.gov/ghgemissions/sources-greenhouse-gas-emissions}}, in 2020 the transportation sector produced 27\% of the total GHG 
emissions in the US, mostly generated from burning fossil fuels by cars, trucks, ships, trains, and planes. Domestic statistics issued by the UK 
government\footnote{\url{https://www.gov.uk/government/statistics/transport-and-environment-statistics-autumn-2021/transport-and-environment-statistics-autumn-2021}} 
confirm that the transportation sector generated 27\% of the total GHG emission. The majority (91\%) came from road transport vehicles, where the biggest 
contributors were cars and taxis. Furthermore, data provided by the European Environment 
Agency\footnote{\url{https://www.eea.europa.eu/data-and-maps/data/data-viewers/eea-greenhouse-gas-projections-data-viewer}} highlight that in the EU more 
than 22\% of the GHG emissions came from the transportation sector.

Despite technical advances have made available a range of options for sustainable mobility, there are still important obstacles that must be overcome for 
their mass adoption. Among such options, Electric Vehicles (EVs) are considered one of the major directions to reduce the environmental impact of people 
mobility and make urban areas more sustainable. In the 2021 edition of the Global EV Outlook 
2021\footnote{\url{https://www.iea.org/reports/global-ev-outlook-2021}}, the International Energy Agency pointed out that at the end of 2020 the global EVs 
stock hit 10 millions units, with 3 millions newly registered EVs. Europe was the fastest growing market, with a sales share equal to 10\% and some leading 
countries, such as Norway, which registered a record high sales share of 75\%. This trend was accelerated by many countries of the European Union through 
substantial financial incentives. However, the decision of potential EV buyers is still strongly affected by two major issues. On one hand, the  purchase 
cost of an EV is still higher than that of a traditional internal combustion engine vehicle. On the other hand, the limited travel range of an EV and the 
long charging time are well-known to generate anxiety in the potential buyers \citep[e.g.,][]{PevBabcarGhiKetPod20}. In fact, the willingness of 
drivers to purchase an EV strongly depends on the availability of charging stations nearby their points of interests (e.g., home and work). As the number 
of charging stations is growing, thanks to public and private investments, the location problem of such stations has attracted much attention  (see 
Section \ref{sec:lit}).

There are a number of factors that make the location of charging stations substantially different from other, more classical, location problems, in 
particular the choice of the charger to install (e.g., slow, quick, fast), and the characteristics of the charging demand.

The type of charger is a key factor to be taken into account, as it impacts the charging time. As of the end of 2021, there exist three main types of 
charger \citep[see][]{Mol21}. Level 1 chargers, also referred to as \textit{slow} chargers, use common 120-volt outlets, and can take up to 40 hours to raise the 
level  of a standard battery EV (with a 60 kWh sized battery) from 10\% to 80\% of the capacity.
These chargers are most suitable for private usage.  Level 2 chargers, sometimes called \textit{quick} chargers, can charge up to 10 times faster than a 
level 1 charger, and are the most commonly used types for daily EV charging \citep[see][]{Mol21}. Given the same battery characteristics mentioned above, 
the charging time is about 4.5 hours. The level 3 or \textit{fast} chargers can reduce the charging time to  40 minutes or even less.
For a comprehensive study regarding the state of the art on charging stations, the interested reader can refer to \citet{pareek2020electric}. The type of 
charger demanded by EVs is affected by the urban layout.  For example, slow chargers will be demanded in residential areas so that EVs can be recharged 
over the night at low cost \citep[an interesting study of the factors influencing the charging demand is provided in][]{wolbertus2018fully}.

In the classical location models a customer is characterized  by the distance from any potential location and by a single quantity - a measure of the 
demand. The models do not consider a temporal dimension of the problem which basically corresponds to assuming that the demand is uniformly distributed 
over the time period of interest of the location decision. On the contrary, the charging demand of EVs fluctuates over time, with peaks of demand in 
periods of time where the traffic volume is high.  Neglecting the demand dynamics may lead to solutions where the charging capacity deployed is not 
sufficient to satisfy the demand during the peak times. 

In this paper, we study the problem of determining an optimal deployment of charging stations for EVs within an urban environment. 
Different types of 
chargers have to be located in pre-defined potential locations, modeled as nodes of a network. 
The urban area is partitioned in sections.
A customer is associated with each section of the urban area. 
Its demand in a certain time 
interval is the number of EVs in that section that need to be recharged.
The customer is located in the center of gravity of the section and is modeled as a node of the network. 
The urban area is also partitioned in zones (e.g., commercial, industrial, or residential) 
which have different needs in terms of minimum  number of each 
type of charger deployed in the zone.

We have to determine, for each type of charger and each potential location, the number of chargers to be deployed. 
Two criteria have a key role in this 
location problem: the cost of installing the chargers and the distance the customers have to travel to be recharged.

We present, over a discretized time horizon, an optimization model that introduces a temporal dimension which, 
to the best of our knowledge, has never been introduced in the literature on location problems and captures the dynamics of the charging 
demand. Assuming that a charger can take more than one period to fully recharge an EV, the proposed multi-period formulation includes 
constraints to keep 
track of the usage of chargers across consecutive time periods and to ensure that no other vehicles are assigned to any occupied charger. 
This novel 
approach guarantees a correct sizing of the solution, in terms of number of stations opened and number of chargers installed,
 and ensures that the demand is completely satisfied in all time periods.
In order to assess the value of introducing the temporal dimension in the location problem, 
which makes the optimization model more complex, we present a single-period 
optimization model that captures the same specificities of the problem but ignores the temporal aspect. 
In both models, the objective is the minimization of a convex combination of two terms: the total cost of deploying the charging stations 
and installing the chargers, and the average distance traveled by the customers to reach the assigned 
charging station.
The two optimization models turn out to be Mixed 
Integer Linear Programming (MILP) problems.
We compare the two models through a theoretical and a computational analysis. We show, through worst-case analysis, that a solution to the 
single-period model may fail to satisfy a large portion of the charging demand. 
Extensive computational experiments are run on different classes of 
randomly generated instances. The results
confirm the importance of 
explicitly considering the dependence on time of the demand. 
In fact, the single-period model is based on the common assumption that the charging demand is 
uniformly distributed across the planning horizon. 
In an application context such as the one at hand, where the demand fluctuates significantly during the 
day and across different zones of the same urban area, 
the single-period model produces solutions that are not capable of serving a large portion of the 
charging demand, especially in those time periods where the demand is prominently concentrated. 
The computational experiments also include a 
parametric analysis of the relative weight assigned to the objective function components.

\textbf{Structure of the paper.} The remainder of the paper is organized as follows. In Section \ref{sec:lit}, 
the literature most closely related to our 
research is reviewed and the contribution of this paper is highlighted. 
In Section \ref{sec:Models}, after the presentation of the single-period extension 
of a classical location model, we provide the multi-period mathematical formulation.
In Section \ref{sec:WCAnalysis},
 we analyze the worst-case performance of the single-period model in terms of  portion of unsatisfied charging demand.
 Section \ref{sec:res} reports extensive computational 
experiments conducted on instances generated to resemble demand dynamics frequently observed in different zones of a city. 
Finally, some concluding remarks are outlined in Section \ref{sec:conc}.

\section{Literature review} \label{sec:lit}

The problem of determining an optimal location and size of charging stations for EVs has recently attracted an increasing academic attention. Recent 
overviews of the main modeling and algorithmic approaches employed in this research area are available in \citet{DebTamKalMah18}, \citet{ZhaLiuZhaGu19}, 
and \citet{Kch21}. For a general introduction on location problems the interested reader can refer to \citet{laporte2019introduction}. In the following, we focus on the papers that are most closely related to our research, and refer the interested reader to the above-mentioned surveys and the references cited therein.

A first broad classification of the literature is based on the type of network considered \citep[cf.][]{DebTamKalMah18}. When only the \emph{distribution 
network} is considered, the optimal location of charging stations must consider the potential adverse effects on the power grid, as an inappropriate 
placement of charging stations can be a threat to the power system security and reliability. On the other hand, when only the \emph{transportation network} 
is taken into account, the main issue is to determine an optimal location of charging stations over a road network. This paper lies in the latter category. 
Within this category, the related literature can be further classified into two main streams of models
called flow-based and node-based demand models \citep[e.g., see][]{Kch21}. 
In the literature, the majority of the research efforts are devoted to the flow-based demand models, whereas the number of papers adopting a node-based 
approach is still relatively limited. To the best of our knowledge, \citet{anjos2020increasing} are the only authors that integrated, within the same 
optimization model, both a node-based and a flow-based approach. 
The flow-based demand models are best suited for modeling long-haul (e.g., inter-urban) journeys where accounting for the limited driving 
range of EVs is important \citep[cf.][]{anjos2020increasing}. Contributions to this line of research can be found, for example, in \citet{kuby2005flow}, \citet{mirhassani2013flexible}, 
\citet{yildiz2016branch}, and \citet{HosMirHoo17}. 
The present paper adopts a \emph{node-based} demand model.

In the class of node-based demand models, drivers demanding to charge their EVs are associated with one/few fixed locations, which represent, for 
instance, their residence, workplace or specific service facilities (such as commercial activities). This approach is best suited for urban  
settings. In fact, in such case EVs do not move much from the location where they need to be charged and their limited driving range  can be neglected 
\citep[cf.][]{anjos2020increasing}.
The most common modeling approaches applied in the literature are based on the extension of classic discrete location models (e.g., location-allocation as 
in \citet{zhu2016charging}, set covering as in \citet{HuaKanZha16}, and maximum coverage problems as in \citet{DonMaWeiHay19}) to incorporate technical 
constraints specific to EVs.

Characteristics of the charging demand (such as the population size, the penetration rate of EVs, the type of zone, and the time of the 
day) are known to have a crucial impact on the optimal location of charging stations. To position the present paper within the literature, we
 classify the mathematical formulations into \emph{single-period}  and \emph{multi-period}. In single-period 
optimization models all the decision variables are time independent. 
Although the spatial-temporal distribution of the charging demands is described by different authors 
\citep[e.g., see][and the references cited therein]{YiZhaLinLiu20}, only few authors have proposed multi-period optimization models 
where the allocation of the demand to the charging 
stations is time-dependent. 
The related stream of literature can be classified according to the length of 
the planning horizon considered. A long planning horizon is considered by some authors. The 
basic rationale of these  models is that locating charging stations is a long-term strategic decision. As a consequence, during these long 
periods of time the technology available, as well as the charging demand, may change significantly. Along this line of research, we mention the 
paper by \citet{anjos2020increasing} where it is assumed that the locating decisions taken in a period have an impact on the charging demand in the subsequent 
periods. In fact, potential EV buyers are influenced by the  availability of charging 
opportunities.
Some papers have proposed  multi-period optimization models that consider a short horizon, 
usually a day,  divided in  time periods, usually hours. 
 Our research belongs to this category of papers. 
 
 To the best of our knowledge, \citet{cavadas2015mip} are the first authors to recognize the importance of incorporating into an optimization model the 
dynamics of the charging demand across the day. The aim of the proposed multi-period model is 
the maximization of the total demand served, subject to a constraint on the budget available. The authors consider only one type of charger (i.e., a slow 
type) and the sizing of the charging stations is not part of the optimization. 
In the model we present in this paper, we address these shortcomings by considering multiple types of 
chargers and optimizing the quantities installed in each opened station.
\citet{RajSad17} estimate the charging demand of EVs in different zones of a city and at different hours. The authors consider the deployment of an unlimited number of fast chargers only and propose a  non-linear optimization model that includes three cost components: the total opening cost, the total cost for  the drivers to reach the assigned charging stations, and the cost of connecting the charging stations to the electric grid substations.
The variability of the demand across the day is taken into consideration when determining the number of chargers to install. Nevertheless, the variables assigning EVs to stations are not time-dependent, and, hence, drivers demanding to charge their EVs at different hours are all assigned to the same station.
In our paper, we allow the demand arising from the same location during the day to be assigned to different stations, depending on the evolution of the overall demand and the available cherging resources. Moreover, we consider different types of chargers.
Both short-term and long-term  decisions are considered in \citet{QudKabMar19}. The main long-term decisions are  related to the year, the location, and the type of charging stations to open. The short-term decisions are mainly related to the amount of power (provided by different sources, such as electric grid and renewable sources) to satisfy the hourly charging demand at a given location. Compared to our research, the drivers are, indirectly, pre-assigned to a charging station and, hence, the assignment is not part of the optimization model. The authors cast the problem as a two-stage stochastic programming model.
\citet{LiJen22} present an optimization model based on the concept of charging opportunities, which is measured through the time an individual stays at a given location within a day. The authors separate the charging opportunities into home and non-home (i.e., public) categories, and allow the same individual to charge 
the EV multiple times at different locations. The proposed optimization model  determines the number of home and non-home chargers to install, as well as 
the times and locations for each individual to charge the EV. The model aims at minimizing the sum of the annual electricity cost for charging the EVs and
 the total cost of locating the home and non-home chargers. The number of chargers that can be installed in each location (called region by the 
authors) is unlimited. 

Finally, we mention the growing body of literature that addresses the problem of determining an optimal location of charging stations for EVs 
in car-sharing systems \citep[e.g., cf.][]{BraKahLei17, brandstatter2020location, BekBoyZog21}. Although such problem has some characteristics in common 
with ours, it includes some operational characteristics that make it considerably different, for example the  decisions 
about the number of EVs to acquire, the relocation of the EVs among stations, and the assumption that charging occurs only between two consecutive 
trips.

\noindent \textbf{Contributions of the paper.}
The contributions of this paper to the literature can be summarized as follows.
\begin{itemize}
\item[\checkmark] We present a node-based multi-period optimization model for the location of charging stations that captures the dependence on time of the charging demand;
\item[\checkmark] the multi-period model takes into account  several characteristics of the real problem: 
multiple types of chargers (each with its own charging speed and installation cost),
 the capacitated nature of the charging stations (in terms of maximum number of chargers that can be installed),
a minimum number of chargers to be installed in different zones (e.g., commercial, residential, industrial);
\item[\checkmark]  the multi-period model is compared to a single-period model through a worst-case analysis;
\item[\checkmark] extensive  computational experiments are presented that show, in particular, the importance of incorporating the dependence on time of the charging demand.
\end{itemize}

\section{Problem definition and mathematical formulations}\label{sec:Models}

In this section, we first provide a general description of the location problem along with the  notation that is common to the two optimization models that 
will follow. Then, the  single-period MILP model is presented, together with the notation that is specific for the model, followed by the multi-period 
formulation.

We consider the problem of determining, in an urban area, an optimal location of charging stations for EVs, along with the type and number of chargers to 
deploy in each station.  A maximum number of chargers, of each type and in total, can be deployed in each station. The location for any station can be 
selected from a pre-defined set of potential locations. We introduce a complete bipartite network $G=(\mathcal{I} \cup \mathcal{J},A)$, where 
$\mathcal{I}=\{1, 2, \dots, I\}$ is the set of demand nodes and $\mathcal{J} = \{1, 2, \dots, J\}$ is the set of potential locations for the stations. Let 
$c_{ij}$ be the travel distance from demand node $i$ to station $j$.

A fixed opening cost $F_j$ is associated with each station $j$. The opening cost does not include the cost of the chargers.  We denote as $\mathcal{K} = 
\{1, 2, \dots, K\}$ the set of types of chargers considered, and as $f_{jk}$ the cost of installing one charger of type $k \in \mathcal{K}$ in location $j 
\in \mathcal{J}$.
Let $u_{jk}$ be the maximum number of chargers of type $k$ that can be installed in station $j$. Similarly, $u_j$ denotes the maximum number of chargers 
that can be installed in total in station $j$. The latter two parameters define, implicitly, the maximum charging capacity of station $j$.

Each node $i$ is the center of gravity  of a section of the urban area where the demand of the section is measured as the number of EVs that need to be 
recharged. We will introduce later, for each of the two optimization models, the planning horizon and  the notation for the demand of a customer. For the 
sake of brevity, hereafter we refer to each potential location $j$ simply as station $j$. The  demand must be entirely satisfied by the chargers that will 
be deployed.

To take into account that different parts of the urban area have different needs in terms of type of charger desired, the urban area is  partitioned in 
zones (e.g., commercial, residential, industrial). We denote by $\mathcal{L} = \{1, 2, \dots, L\}$ the set of zones. We assume that, based on some 
preliminary analysis, in each zone $\ell \in \mathcal{L}$ a minimum percentage $\rho_{\ell k}$ of chargers of type $k$ must be deployed. Each station $j\in 
\mathcal{J}$  belongs to a zone as well as each customer $i \in \mathcal{I}$. Thus, the  zones imply a partition of both the stations and the demand points. This partition does not restrict the allocation of demand to stations, i.e., a demand point located in a zone can be assigned to a station located in a different zone.

Two criteria have a key role in this location problem: the cost of opening the stations and installing the chargers and the distance the customers have to 
travel to be recharged. The objective function we consider, to be minimized, is a convex combination of these two criteria.
The optimization problem is aimed at determining, for each type of charger and each station, the number of chargers to be deployed in such a way that the 
objective function is minimized.

Both MILP models include the following decision variables. Let $z_j \in \{0, 1\}$, with $j \in \mathcal{J}$, be a binary variable that takes value 1 if 
station $j$ is opened, and 0 otherwise. Let $y_{jk} \in \mathbb{Z}_+$, with $j \in \mathcal{J}$ and $k \in \mathcal{K}$, be an integer variable that 
represents the number of chargers of type $k$ installed in station $j$.

\subsection{A  single-period location model} \label{sec:SP-CFLP}

This section presents a single-period model for the location of the charging stations. The MILP formulation, denoted as \spmodel, is an extension of a 
classical CFL model.
Hereafter, we introduce the  notation needed for the formulation, in addition to the one  introduced above.

We consider a single planning period of length $H$ and denote as $d_i$ the total demand in $i \in \mathcal{I}$, that is, the total number of EVs demanding 
to be recharged in $i$ during $H$. Let $p_k$ denote the average number of EVs fully recharged by one charger of type $k$ during time period $H$. For the 
sake of simplicity, we assume that $p_k$ does not depend on the type of EV.

The \spmodel\ model also makes use of the following decision variables. Let $x_{ijk} \in [0, 1]$, with $i \in \mathcal{I}$, $j \in \mathcal{J}$, and $k \in 
\mathcal{K}$, be the fraction of the  demand of node $i$ assigned to a charger of type $k$ in station $j$. Then, the \spmodel\ model can be stated as the 
following MILP:

\begin{equation} \label{M4:OF}
\mbox{[\spmodel]} \quad \mbox{min} \quad  \lambda \cdot \left(\frac{1}{\sum\limits_{i \in \mathcal{I}} d_{i}} \sum_{i \in \mathcal{I}} d_{i} \sum_{j \in 
\mathcal{J}} c_{ij} \sum_{k \in \mathcal{K}} x_{ijk}\right) + (1 -\lambda) \cdot \left(\sum_{j \in \mathcal{J}} F_j z_j + \sum_{j \in \mathcal{J}} \sum_{k 
\in \mathcal{K}} f_{jk} y_{jk}\right)
\end{equation}

\begin{equation} \label{M4:C1}
\mbox{s.t.} \quad y_{jk} \leq u_{jk} z_j \quad j \in \mathcal{J}, k \in \mathcal{K}
\end{equation}

\begin{equation} \label{M4:C2}
\sum_{k \in \mathcal{K}} y_{jk} \leq u_{j} z_j \quad j \in \mathcal{J}
\end{equation}

\begin{equation} \label{M4:C3}
\sum_{j \in \mathcal{J}} \sum_{k \in \mathcal{K}} x_{ijk} = 1 \quad i \in \mathcal{I}
\end{equation}

\begin{equation} \label{M4:C4}
\sum_{i \in \mathcal{I}} d_{i} x_{ijk} \leq p_k y_{jk} \quad j \in \mathcal{J}, k \in \mathcal{K}
\end{equation}

\begin{equation} \label{M4:C5}
x_{ijk} \leq y_{jk} \quad i \in \mathcal{I}, j \in \mathcal{J}, k \in \mathcal{K}
\end{equation}

\begin{equation} \label{M4:C6}
\sum_{j \in A_{\ell}} y_{jk} \geq \rho_{\ell k} \sum_{j \in A_{\ell}} \sum_{k \in \mathcal{K}} y_{jk} \quad k \in \mathcal{K}, \ell \in  \mathcal{L}
\end{equation}

\begin{equation} \label{M4:C7}
z_j \in \{0, 1\} \quad j \in \mathcal{J}; \quad y_{jk} \in \mathbb{Z}_+ \quad j \in \mathcal{J}, k \in \mathcal{K}; \quad x_{ijk} \in [0, 1] \quad i \in 
\mathcal{I}, j \in \mathcal{J}, k \in \mathcal{K}.
\end{equation}


The objective function in \eqref{M4:OF} comprises two terms. The first one represents the average distance traveled by the EVs to reach the assigned 
station. The second term is the total cost of opening the stations and installing the chargers. The two  terms represent criteria of a substantially 
different nature: the first measures the quality of the service provided by the deployed stations and chargers to the drivers, whereas the second the cost 
of the service. The two criteria are weighted  by the trade-off parameter $\lambda \in [0, 1]$, which is used to balance their importance.

Constraints \eqref{M4:C1} and \eqref{M4:C2} limit the number of chargers that can be installed in station $j$. The former set bounds the number of chargers 
of type $k$ to be lower than or equal to $u_{jk}$, whereas the second set of constraints bounds the total number of chargers to be lower than or equal to 
$u_{j}$. Both sets of constraints \eqref{M4:C1} and \eqref{M4:C2} impose that no charger can be installed if station $j$ is not open (i.e., $z_j = 0$). 
Constraints \eqref{M4:C3} ensure that the demand of each node $i \in \mathcal{I}$ is entirely satisfied. Constraints \eqref{M4:C4} guarantee that the 
number of EVs assigned to the chargers of type $k$ deployed in station $j$ is not greater than the charging capacity available (i.e., $p_k y_{jk}$). They 
also impose that no EV can be assigned to a type $k$ of chargers in station $j$ if no charger of that type is available (i.e., $y_{jk} = 0$). Inequalities 
\eqref{M4:C5}, which are redundant in this formulation, are well-known to yield a  tighter Linear Programming (LP) relaxation than the equivalent 
formulation without  them \citep[e.g., see][]{FilGuaSpe21}.
Constraints \eqref{M4:C6} guarantee that the number of chargers of type $k$ installed in zone $\ell$ is at least equal to the minimum percentage 
$\rho_{\ell k}$. Finally, constraints \eqref{M4:C7} define the domain of the decision variables.

\subsection{A multi-period  location model}\label{sec:MP-CFLP}

This section presents the MILP formulation for the multi-period model, henceforth denoted as the \mpmodel\ model, for the problem defined at 
the beginning of this section.

The planning period $H$ of the single-period model is here partitioned into a number $T$ of time periods. For example, if $H$ is a day, we may partition the day in hours. Let $\mathcal{T} = \{1, 2, \dots, T\}$ denote the set of time periods. We denote as $R_k$ the number of consecutive time periods needed 
to completely recharge a car using a charger of type $k$. Note that, similar to $p_k$ for the \spmodel\ model, $R_k$ does not depend on the type of EV but 
only on the type of charger. Furthermore, parameters $p_k$ and $R_k$ are strictly related, as the latter is determined by dividing the length of the time 
horizon by $p_k$, i.e. $R_k = \frac{T}{p_k}$.

The demand of each node $i \in \mathcal{I}$ is no longer identified by a single value ($d_i$ in the \spmodel\ model) but by a time-dependent profile. Let 
$d_i^t$ denote the demand of node $i \in \mathcal{I}$  at the beginning of time period $t \in \mathcal{T}$. A more detailed discussion about the demand 
profiles can be found in Section \ref{sec:IG}. We assume that the demand of a  time period $t$ must be served in that time period, i.e., it cannot be 
postponed to a later time.  We say that a node is served by a charger of type $k$ at time $t$ if a charger is available at time $t$ to start the charging 
which will occupy the charger for a total of $R_k$ time periods. The capacity installed in each station must be sufficient to serve the charging demand 
assigned to that station in a time period and the demand assigned to the station in a previous time period that has not yet completed the charging. 
Finally, let $x_{ijk}^t \in [0, 1]$, with $i \in \mathcal{I}$, $j \in \mathcal{J}$, $k \in \mathcal{K}$, and $t \in \mathcal{T}$, be the fraction of the 
charging demand of  node $i$ to be served at time  $t$ that is assigned to a charger of type $k$ in station $j$.

The \mpmodel\ model is formulated as follows:

\begin{center}
[\mpmodel]
\end{center}
\begin{equation} \label{M1:OF}
\quad  \mbox{min} \quad  \lambda \cdot \left(\frac{1}{\sum\limits_{t \in \mathcal{T}} \sum\limits_{i \in \mathcal{I}} d_{i}^{t}} \sum\limits_{t \in 
\mathcal{T}} \sum\limits_{i \in \mathcal{I}} d_{i}^{t} \sum\limits_{j \in \mathcal{J}} c_{ij} \sum\limits_{k \in \mathcal{K}} x_{ijk}^{t}\right) + (1 
-\lambda) \cdot \left(\sum\limits_{j \in \mathcal{J}} F_j z_j + \sum\limits_{j \in \mathcal{J}} \sum\limits_{k \in \mathcal{K}} f_{jk} y_{jk}
\right)
\end{equation}

\begin{center}
s.t. \eqref{M4:C1}, \eqref{M4:C2}, and \eqref{M4:C6}
\end{center}

\begin{equation} \label{M1:C1}
\sum_{j \in \mathcal{J}} \sum_{k \in \mathcal{K}} x_{ijk}^t = 1 \quad i \in \mathcal{I}, t \in \mathcal{T}
\end{equation}



\begin{equation} \label{M1:C4}
x_{ijk}^t \leq y_{jk} \quad i \in \mathcal{I}, j \in \mathcal{J}, k \in \mathcal{K}, t \in \mathcal{T}
\end{equation}

\begin{equation} \label{M1:C5}
\sum_{i \in \mathcal{I}} \sum_{\tau=0}^{t-1} d_{i}^{t-\tau} x_{ijk}^{t-\tau} \leq y_{jk} \quad j \in \mathcal{J}, k \in \mathcal{K},
t \in \mathcal{T}: t < R_k
\end{equation}

\begin{equation} \label{M1:C6}
\sum_{i \in \mathcal{I}} \sum_{\tau=0}^{R_k - 1} d_{i}^{t-\tau} x_{ijk}^{t-\tau} \leq y_{jk} \quad j \in \mathcal{J}, k \in \mathcal{K}, t \in \mathcal{T}: 
t \geq R_k
\end{equation}


\begin{equation} \label{M1:C8}
z_j \in \{0, 1\} \quad j \in \mathcal{J}; \quad y_{jk} \in \mathbb{Z}_+ \quad j \in \mathcal{J}, k \in \mathcal{K}; \quad x_{ijk}^t \in [0, 1] \quad i \in 
\mathcal{I}, j \in \mathcal{J}, k \in \mathcal{K}, t \in \mathcal{T}.
\end{equation}

The objective function in \eqref{M1:OF} is the multi-period extension of function \eqref{M4:OF}. For each node $i \in \mathcal{I}$, constraints 
\eqref{M1:C1} ensure that the charging demand arising in each time period $t$ is fully satisfied. Akin to the objective function, also inequalities 
\eqref{M1:C4} are the multi-period extension of constraints \eqref{M4:C5}.

Constraints \eqref{M1:C5} and \eqref{M1:C6} guarantee that the number of EVs that are charging in time period $t$ at a charger of type $k$ in station $j$ 
is smaller than or equal to the number of available chargers of that type (i.e., $y_{jk}$). Note that the second sum in \eqref{M1:C5} and \eqref{M1:C6} is 
used to keep track of the EVs that started to recharge in a previous time period but have not completed the charging in $t$. Constraints \eqref{M1:C5} 
are defined for the first time periods in the planning horizon (such that $t < R_k$), whereas \eqref{M1:C6} are defined for the remaining time periods. 
Finally, constraints \eqref{M1:C8} define the domain of the decision variables.

\section{Worst-case analysis}\label{sec:WCAnalysis}

In this section, we analyze the worst-case performance of the \spmodel\ model in terms of the demand that cannot be satisfied if the optimal solution 
produced is implemented in a context where the demand fluctuates  over time. In fact, in this case if an optimal solution to the \spmodel\ model is 
implemented, there is no guarantee that all the charging demand is satisfied. As the \spmodel\ model implicitly assumes that the charging demand is 
uniformly distributed across the planning horizon, when the demand fluctuates over time, there may be peak time periods where the chargers installed are 
not sufficient.

\begin{theorem}\label{prp:errorbound}
When an optimal solution of the \spmodel\ model is implemented, the fraction of the demand that does not find an available charger to be served may be up 
to $1-\frac{1}{T}$, where $T$ is the number of time periods of the planning horizon. This bound is tight.
\end{theorem}
\proof
To prove the theorem, we build the following instance.

Recalling constraint \eqref{M4:C3} and summing up all constraints \eqref{M4:C4}, the following chain of inequalities holds:
\begin{equation}
\sum_{i \in \mathcal{I}} d_i \stackrel{\eqref{M4:C3}}{=} \sum_{i \in \mathcal{I}} d_i \sum_{j \in \mathcal{J}}\sum_{k \in \mathcal{K}} x_{ijk} = \sum_{j 
\in \mathcal{J}} \sum_{k \in \mathcal{K}} \sum_{i \in \mathcal{I}} d_i x_{ijk} \stackrel{\eqref{M4:C4}}{\leq} \sum_{j \in \mathcal{J}} \sum_{k \in 
\mathcal{K}} p_k y_{jk}
\label{eq:total-sum}
\end{equation}%
for any feasible solution to the \spmodel\ model.

Consider an instance where the travel distances are all negligible compared to the fixed opening and installing cost. In this situation, the \spmodel\ 
model would open the minimum number of charging stations and install the minimum number of chargers that are strictly necessary to satisfy the total 
demand. As a consequence, the value of the right-hand side of the rightmost inequality in \eqref{eq:total-sum} would be as small as possible.

Additionally, suppose there is a single type of charger (i.e., $K=1$) and that the total demand $\sum_{i \in \mathcal{I}} d_i$ is a multiple of $p_1$. 
Recall that the latter parameter represents the  number of EVs fully recharged by one charger during the planning horizon. Note that it can be determined 
by dividing the number of time periods $T$ by the  number of consecutive time periods needed to completely recharge an EV (i.e., $R_1$). Hence, at 
optimality, the inequality in \eqref{eq:total-sum} can be reformulated as follows:
\begin{equation}
\sum_{i \in \mathcal{I}} d_i = \sum_{j \in \mathcal{J}} p_1 y_{j1} = \frac{T}{R_1} \sum_{j \in \mathcal{J}} y_{j1}.
\end{equation}%
Thus, the total number of chargers deployed is $\sum_{j \in \mathcal{J}} y_{j1} = \frac{R_1 \sum_{i \in \mathcal{I}} d_i}{T}$, which, assuming that $R_1 = 
1$, becomes $\sum_{j \in \mathcal{J}} y_{j1} = \frac{\sum_{i \in \mathcal{I}} d_i}{T}$.

Consider an extreme situation where the whole demand $\sum_{i \in \mathcal{I}} d_i$ arises in one time period, say $\hat t$, whereas it is zero in the 
remaining periods. The demand that can be satisfied in such time period is equal to the number of chargers installed (i.e., $\sum_{j \in \mathcal{J}} 
y_{j1}$). Given the assumptions above, this value is also equal to $\frac{\sum_{i \in \mathcal{I}} d_i}{T}$. Thus, the amount of demand that does not find 
an available charger is equal to:
\[
\sum\limits_{i \in \mathcal{I}} d_i - \frac{\sum\limits_{i \in \mathcal{I}} d_i}{T}.
\]

\begin{figure} [H]
\centering
		\includegraphics[width=0.9\linewidth, trim={.15cm 0 0 0}, clip=true]{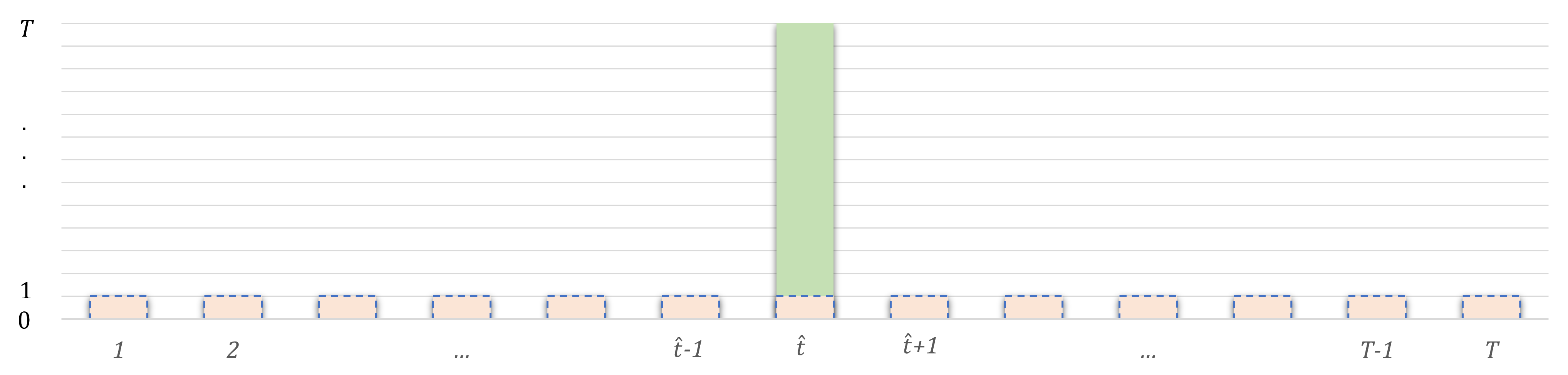}%
	\caption{An instance where the  demand arises in time period $\hat t$ (green bar). The pink bars show a uniform distribution of the demand across 
the planning horizon, as implicitly assumed by the \spmodel\ model.}%
	\label{fig:WC_bound}%
\end{figure}

The statement follows.
\prend
Figure \ref{fig:WC_bound} illustrates the construction for the special case where $\sum_{i \in \mathcal{I}} d_i = T$, which implies that $\sum_{j \in 
\mathcal{J}} y_{j1} = 1$.
The whole demand, equal to $T$, arises in time period $\hat t$ (green bar), whereas it is zero in the remaining periods. The \spmodel\ model assumes that 
such demand is uniformly distributed across the planning horizon (pink bars), and hence it opens one station equipped with one charger. As a consequence, 
the charging demand that is not satisfied is $T-1$ or, in percentage, $\frac{T-1}{T}$.

\section{Experimental Analysis} \label{sec:res}

This section is devoted  to the presentation and discussion of the computational experiments. 
They were conducted on a Workstation HP Intel(R)-Xeon(R) at 
3.5GHz with 64 GB RAM (Win 10 Pro, 64 bits). The processor is equipped with 6 physical cores, and all threads were used while solving each instance. The 
MILP models were implemented in Java, compiled within Apache NetBeans 12.3, and solved by means of CPLEX 20.1. Each instance was solved with a CPU time 
limit of 3,600 seconds. All other CPLEX parameters were set at their default values.

The section is organized as follows. First, we present the testing environment we used in our experiments, 
then we compare the optimal solutions for two 
illustrative examples generated according to two different urban structure models, and finally we provide detailed computational results 
comparing the 
solutions produced by the single-period and the multi-period models.

\subsection{Testing environment}\label{sec:Testing}

The generation of the charging demand and potential station locations follows the procedure described in Section~\ref{sec:IG}. All the remaining parameters 
defining the testing environment are detailed in Section~\ref{sec:RemainingParam}.

\subsubsection{Spatial and temporal charging demand generation} \label{sec:IG}

As far as the urban structure is concerned, we considered two classic models, the concentric zone model and the sector model.
The \emph{concentric zone model} was proposed in 1925 by sociologist Ernest Burgess on the base of his human ecology theory, and was initially applied 
to the city of Chicago \citep[cf.][]{2008-Bur}. It is, perhaps, the first theoretical model used to explain urban social structures. The model depicts 
urban land usage as concentric rings: the business district is located in the center, whereas the remainder of the city is expanded in rings, each 
corresponding to a different land usage (such as industrial or residential).
The \emph{sector model} was proposed in 1939 by land economist Homer Hoyt \citep[see][]{1939-Hoy}. It is a modification of the Burgess' model where the 
city zones devoted to a specific land usage (e.g., business, residential, and productive) develop in sectors expanding from the original city center.
Though the actual structure of modern cities can hardly be captured by models as simple as Burgess' and Hoyt's, they are the basis of more complex 
structures \citep{2012-Hall} and, on the other hand, can simplify the interpretation of the results. For these reasons, we considered two classes of 
instances, each associated with one urban model. Hereafter, the two classes are referred to as the \emph{concentric ring} (COR) instances and the 
\emph{sector} (SEC) instances. In both cases, we assume the urban structure comprises three possible zones: \emph{commercial}, \emph{residential}, 
\emph{industrial}. With a little abuse of notation, we denote the set of zones as $\mathcal{L} = \{C, R, I\}$, where $C$, $R$, and $I$ refer to the 
commercial, residential, and industrial zones, respectively. Each zone is characterized by a different pattern of the charging demand during the planning 
horizon, as it will be detailed later. We consider as planning horizon a day, discretized into $T = 24$ hours (i.e., the time periods).

In the COR instances, we assume the commercial zone is the central circle with ray 1000, the residential zone is a ring in the middle around the commercial 
zone with outer ray 2000, and the industrial zone is a most outer ring around the residential zone with outer ray 3000. In the SEC instances, we assume 
that commercial, residential, and industrial zones correspond to three slices of identical size that partition a circle of ray 3000.

\begin{figure} 
	\begin{subfigure}{0.475\textwidth}
		\centering
		\includegraphics[width=0.9\linewidth]{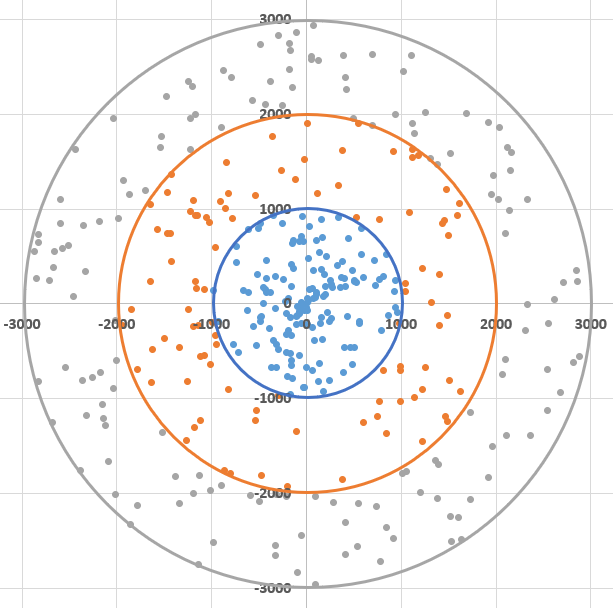}%
		\caption{COR instance \label{subfig:e}}
	\end{subfigure}
	\begin{subfigure}{0.475\textwidth}
		\centering
		\includegraphics[width=0.9\linewidth]{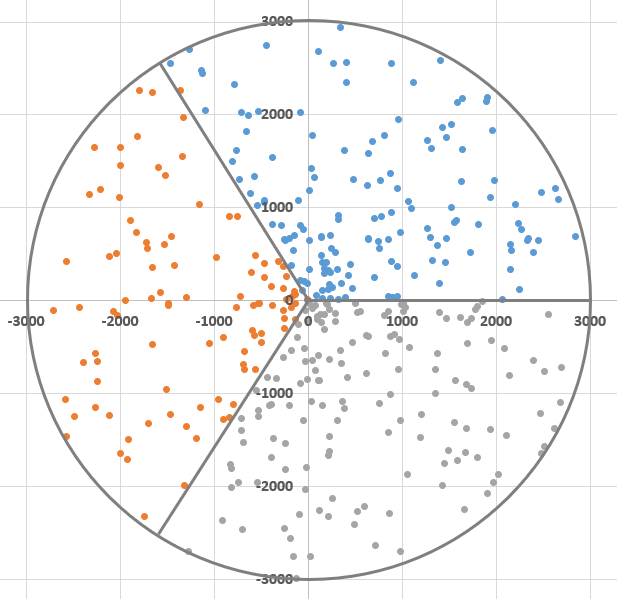}%
		\caption{SEC instance}
	\end{subfigure}
	\caption{Demand nodes generation for a COR (left) and a SEC (right) instance. Blue points are located in the commercial zone, orange points in the 
residential zone, and gray points in the industrial zone.}%
	\label{fig:zones}%
\end{figure}

Then, for each zone, we uniformly generate the same number of demand nodes. More exactly, given the total number of demand nodes to be generated, such 
value is divided by three to obtain, after an integer rounding whenever necessary, the number of demand nodes to generate in each zone. 
Figure~\ref{fig:zones} gives an example.

Both in COR and SEC instances, a given number $J$ of potential stations is uniformly generated over the total area (i.e., the circle with ray 3000).

For each pair of demand node $i$ and potential station $j$, parameter $c_{ij}$ is computed as the Euclidean distance between the two nodes.

Concerning the demand pattern, this is specific for each zone. In the commercial zone we assume there is a high density of offices, shops, pubs, 
restaurants, and hotels. Hence, we expect a high demand with two peaks, the first one at the beginning of the working day and the second one in the 
afternoon, higher than the first peak and slowly decreasing in the evening hours. In the industrial zone, we expect a high demand in the morning with one 
peak around lunch time and an almost null demand during the night. Finally, in the residential zone, we assume the presence of a high density of private 
houses, with a comparatively low demand in the morning and a peak in the evening and early night hours. These assumptions are consistent with several 
studies regarding the spatial-temporal distribution of the charging demands observed in urban areas \citep[see, e.g.,][]{YiZhaLinLiu20, 2021-StrStrGoh}.

\begin{figure}
	\begin{subfigure}{0.9\textwidth}
		\centering
		\includegraphics[width=0.9\linewidth]{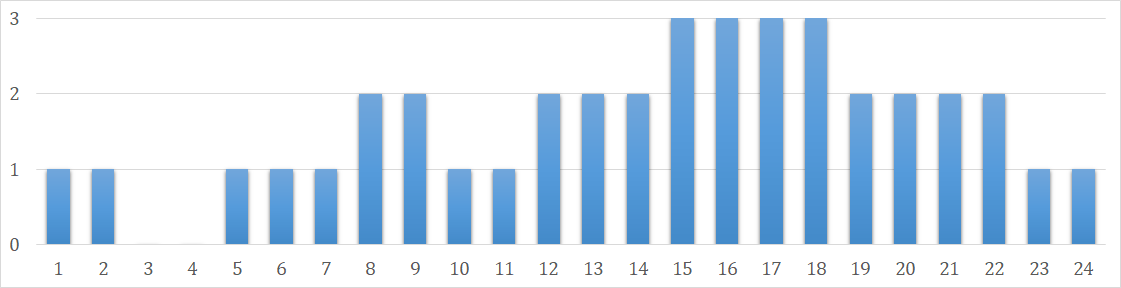}%
		\caption{Commercial}
	\end{subfigure}
\begin{subfigure}{0.9\textwidth}
	\centering
	\includegraphics[width=0.9\linewidth]{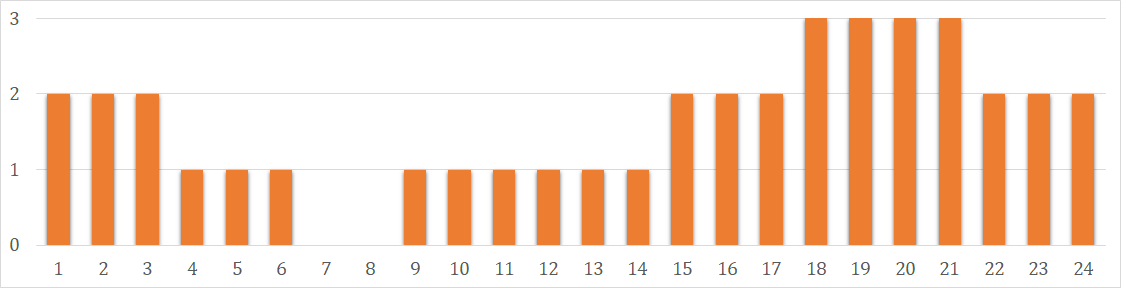}%
	\caption{Residential}
\end{subfigure}
	\begin{subfigure}{0.9\textwidth}
		\centering
		\includegraphics[width=0.9\linewidth]{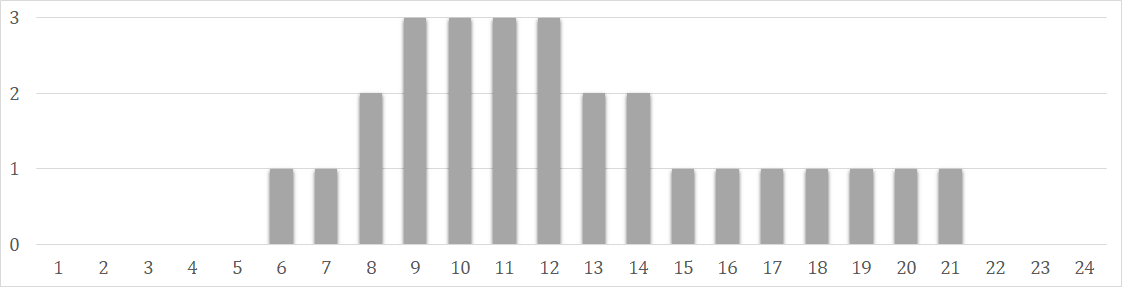}%
		\caption{Industrial}
	\end{subfigure}

	\caption{Basic hourly demand profile in each zone.}%
	\label{fig:profiles}%
\end{figure}

To generate the charging demand according to these patterns, we proceed as follows. We first generate three basic profiles of demand that mimic the 
patterns described above for the commercial, industrial, and residential zones, respectively. Such basic profiles are built according to the four standard 
demand levels shown in Table~\ref{tab:levels}. The resulting basic demand profiles are depicted in Figure~\ref{fig:profiles}.
\begin{table}
	\centering
	\begin{tabular}{r|cccc}
		Level   & 0 & 1 & 2 & 3 \\ \hline
		Demand & null & low & medium & high
	\end{tabular}
	\caption{Standard levels of hourly demand.}
	\label{tab:levels}
\end{table}
For each $i \in \mathcal{I}$ and $t=1, \dots, 24$, we randomly generate an initial demand value $\tilde d_i^t$ from a Poisson distribution with mean (and 
variance) equal to the standard level assigned to $i$ and $t$ in the corresponding basic profile. For example, if $i$ is in the commercial zone and $t=8$, 
then $\tilde d_i^t$ is a realization of a Poisson with mean 2, cf. Figure~\ref{fig:profiles}(a). We then set:
\[
d_i^t = \left[ \frac{\tilde d_i^t}{\sum_{t \in \mathcal{T}}\tilde d_i^t}\cdot 10\right],
\]%
where $[\cdot]$ denotes the nearest integer rounding operator. In this way, we obtain that: (1) the total daily demand from each demand node is around 10; 
(2) the total demand in each zone is consistent with the corresponding basic demand profiles shown in Figure~\ref{fig:profiles}.


\subsubsection{Remaining parameters}\label{sec:RemainingParam}

We generated a set of instances by varying the number of demand nodes $I$, of potential stations $J$, and of the maximum number of chargers to install in each 
station $u_j$. All the remaining parameters take the same value across all instances. The name of each instance is $I\_J\_u_j$, where:

\begin{itemize}
	\item[\checkmark] $I$: The number of demand nodes ranges according to the following values: $I = 50$, 100, 150, 200, 250, and 500.
	\item[\checkmark] $J$: The number of potential stations ranges according to the following values: $J = 10$, 20, 30, 40, and 50.
	\item[\checkmark] $u_j = u$: The maximum number of chargers to install is equal across all the stations, and ranges according to the following 
values: $u_j = 10$, 20, and 30.
\end{itemize}

For example, instance \texttt{50\textunderscore 20\textunderscore 30} comprises 50 demand nodes, 20 potential stations, and parameter $u_j$ is equal to 30. 
Note that the latter parameter is equal for each potential station $j$. We made this choice to simplify the interpretation of the results. For the same 
reason, we decided to set $u_{jk} = u_j = u$, for each $j \in \mathcal{J}$ and $k \in \mathcal{K}$.

Each instance in our testbed has the following common characteristics:

\begin{itemize}
	\item[\checkmark] The planning horizon considered is 1 day, discretized in time periods of one hour length. Consequently, $\mathcal{T} = \{1, 2, 
\dots, 24\}$.
	\item[\checkmark] The cost of opening one charging station is $F_j = $100,000$ \  \forall j \in J$.
	\item[\checkmark] Two types of chargers are considered. Therefore, $\mathcal{K} = \{1, 2\}$, where 1 denotes quick chargers, and 2 stands for fast 
chargers.
	\item[\checkmark] The cost of installing each type of chargers is $f_{j1} = 3,000$  and $f_{j2} = $25,000$ \ \forall j \in 
J$ for quick and fast chargers, respectively. 
	\item[\checkmark] Quick chargers need $R_1 = 4$ hours to fully recharge an EV, whereas fast chargers require $R_2 = 1$ hour. 
Parameter $p_k$ for the \spmodel\ model is determined as: $p_k = \frac{24}{R_k}$.
	\item[\checkmark] The minimum percentage of chargers of each type to deploy in each zone is the following:
	\begin{itemize}
		\item Commercial zone: at least 20\% of quick chargers ($\rho_{C1} = 0.20$), and at least 40\% of fast chargers ($\rho_{C2} = 0.40$).
		\item Residential zone: at least 50\% of quick chargers ($\rho_{R1} = 0.50$), and at least 20\% of fast chargers ($\rho_{R2} = 0.20$).
		\item Industrial zone: at least 25\% of quick chargers ($\rho_{I1} = 0.25$), and at least 25\% of fast chargers ($\rho_{I2} = 0.25$).
		
	\end{itemize}
\end{itemize}

Note that, in both MILP models, the two components in the respective objective function can take very different values, differing even by  orders of 
magnitude. In our experiments we scaled the two components to make them 
comparable in value.\\
We initially considered all possible combinations of the values mentioned above for $I$, $J$, and $u_j$. Subsequently, we ruled out each instance that 
turned out to be infeasible for both MILP models. This situation happened especially for the largest numbers of the demand nodes, the 
smallest numbers of potential stations, and the smallest values of parameter $u_j$. In such cases, the maximum charging capacity, obtained opening all 
potential stations and deploying $u_j$ chargers in each station $j$, turned out not to be sufficient to serve the total charging demand.
All together, we analyzed 67 instances.

\subsection{A comparison between COR and SEC instances}\label{sec:CompareUrbanModels}

\begin{figure}
	\begin{subfigure}{0.5\textwidth}
		\centering
		\includegraphics[width=0.9\linewidth]{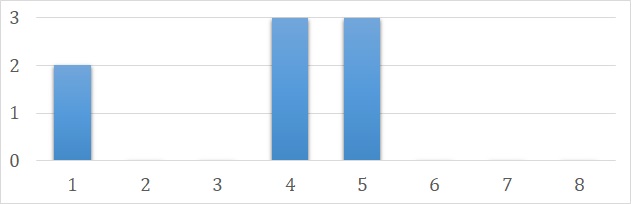}%
		\caption{Commercial}
	\end{subfigure}
	\begin{subfigure}{0.5\textwidth}
		\centering
		\includegraphics[width=0.9\linewidth]{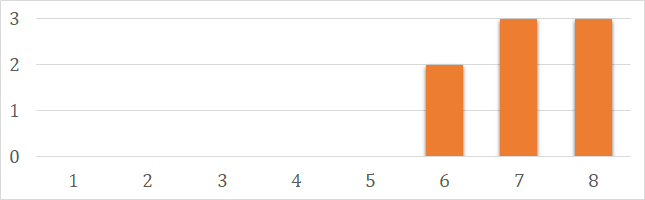}%
		\caption{Residential}
	\end{subfigure}
	\\[0.3cm]
	\begin{center}
		\begin{subfigure}{0.5\textwidth}
		\centering
		\includegraphics[width=0.9\linewidth]{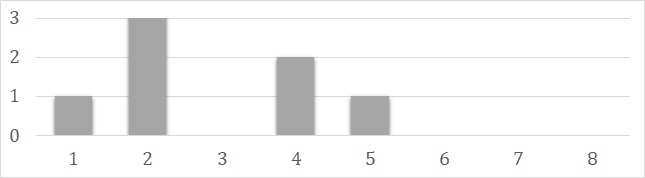}%
		\caption{Industrial}
	\end{subfigure}

\end{center}

	\caption{Illustrative example: Basic hourly demand profile in each zone.}%
	\label{fig:profiles_test}%
\end{figure}

To illustrate the solutions obtained by the two MILP models on the COR and SEC instances, we  discuss the results obtained on two small instances. 
In both instances, the number $I$ of demand nodes is equal to 
21, equally divided among the three zones. The number $J$ of potential stations is 5.
We assumed that the planning horizon comprises 8 time periods, and that the demand profile in each zone is the one depicted in 
Figure~\ref{fig:profiles_test}. These profiles are the same for both the COR and SEC instance.

\begin{sidewaysfigure}

	\begin{subfigure}{0.5\textwidth}
		\centering
		\includegraphics[width=0.8\linewidth]{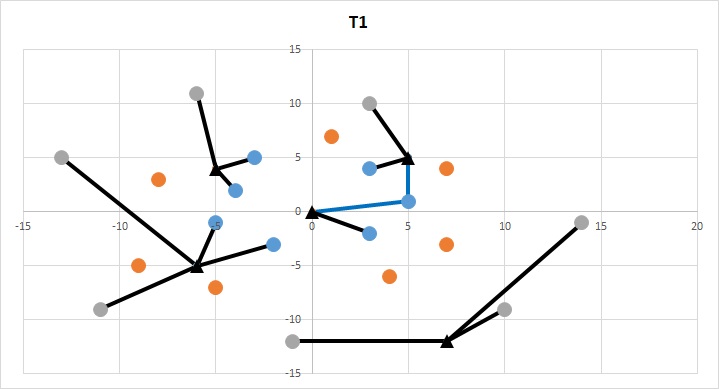}%
		\caption{Time Period 1}
	\end{subfigure}
		\begin{subfigure}{0.5\textwidth}
		\centering
		\includegraphics[width=0.8\linewidth]{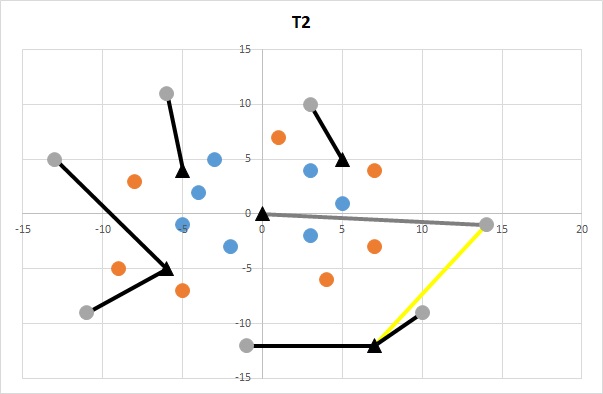}%
		\caption{Time Period 2}
	\end{subfigure}
	\\[0.3cm]
		\begin{subfigure}{0.5\textwidth}
		\centering
		\includegraphics[width=0.8\linewidth]{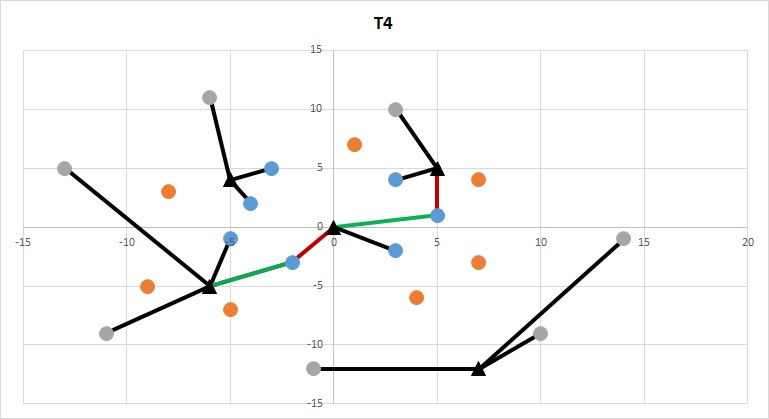}%
		\caption{Time Period 4}
	\end{subfigure}
		\begin{subfigure}{0.5\textwidth}
		\centering
		\includegraphics[width=0.8\linewidth]{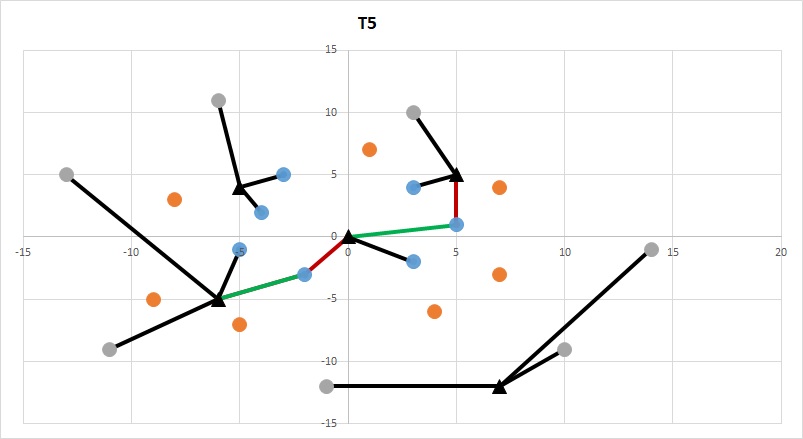}%
		\caption{Time Period 5}
	\end{subfigure}
	\\[0.3cm]
		\begin{subfigure}{0.5\textwidth}
		\centering
		\includegraphics[width=0.8\linewidth]{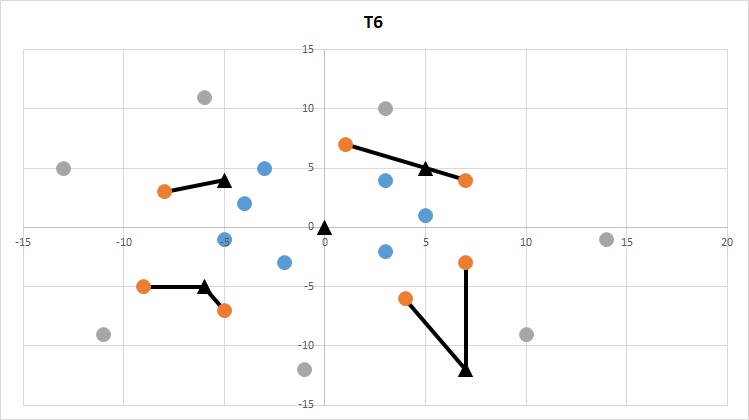}%
		\caption{Time Period 6}
	\end{subfigure}
		\begin{subfigure}{0.5\textwidth}
		\centering
		\includegraphics[width=0.8\linewidth]{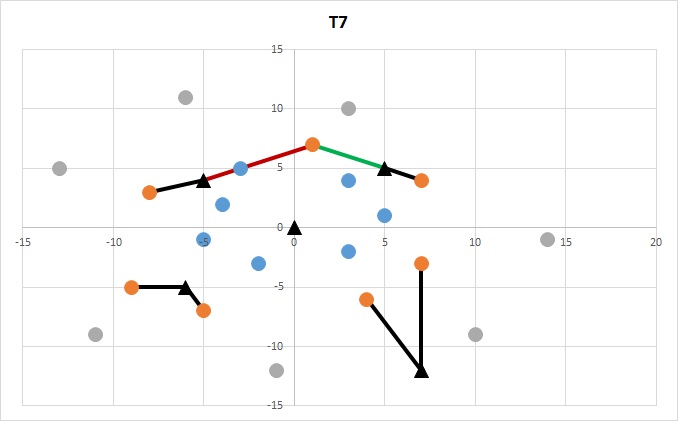}%
		\caption{Time Period 7}
	\end{subfigure}
\caption{COR instance: An optimal solution to the \mpmodel\ model.}
	\label{fig:M1_cerchi}
\end{sidewaysfigure}

\begin{figure}
	\centering
	\includegraphics[scale=0.75]{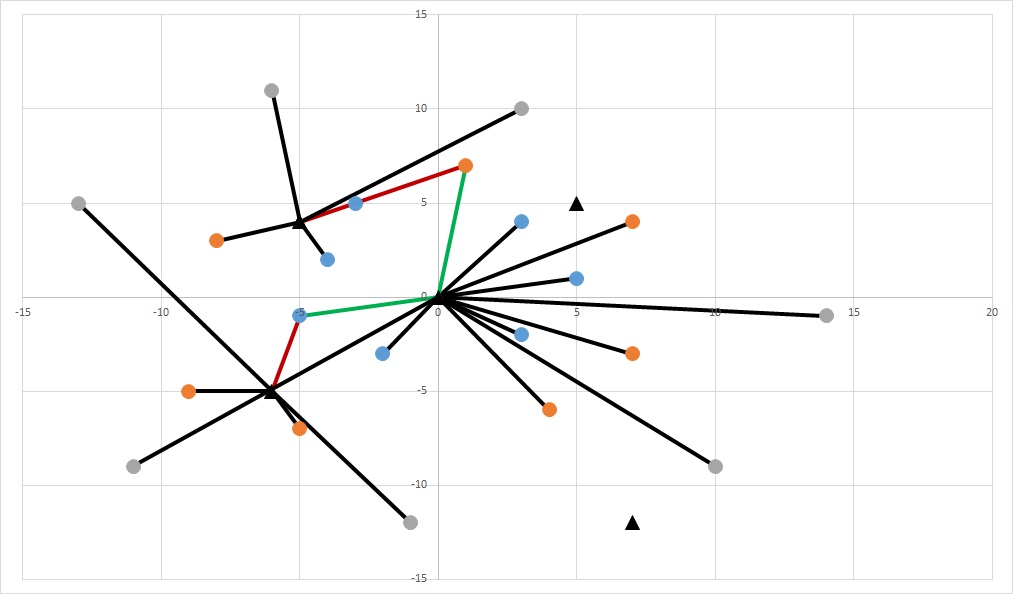}
		\caption{COR instance: An optimal solution to the \spmodel\ model.}
		\label{fig:M4_cerchi}
\end{figure}

An optimal solution produced by the \mpmodel\ model for the COR instance is depicted in Figure~\ref{fig:M1_cerchi}, where demand nodes are colored circles 
(blue for the commercial zone, orange for the residential zone, grey for the industrial zone) and potential locations are black triangles. Moreover, the 
color of the edge connecting a demand node to a black triangle represents the fraction of the charging demand assigned to the station thereby opened (black 
= 100\% of the demand, yellow 75\%, red 66\%, blue 50\%, green 33\% and gray 25\%). Note that Figure~\ref{fig:M1_cerchi} shows the assignments concerning 
only the significant time periods. In other words, the assignments in time periods 3 and 8 are not reported, the former since there is no demand in that 
time period, the latter because it is identical to time period 7. Figure~\ref{fig:M4_cerchi} displays an optimal solution to the \spmodel\ model for the same instance.

\begin{sidewaysfigure}
	\begin{subfigure}{0.5\textwidth}
		\centering
		\includegraphics[width=0.73\linewidth]{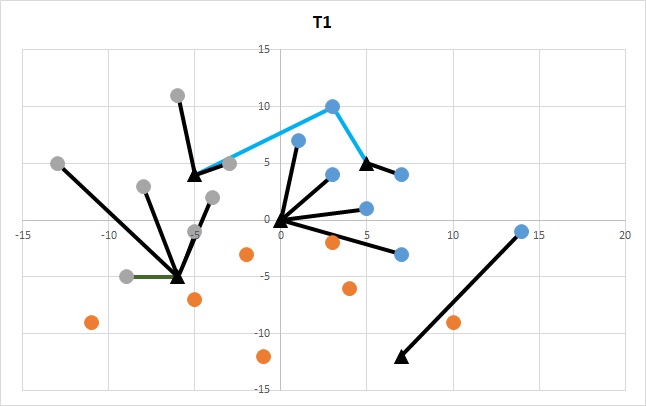}%
		\caption{Time Period 1}
	\end{subfigure}
	\begin{subfigure}{0.5\textwidth}
		\centering
		\includegraphics[width=0.73\linewidth]{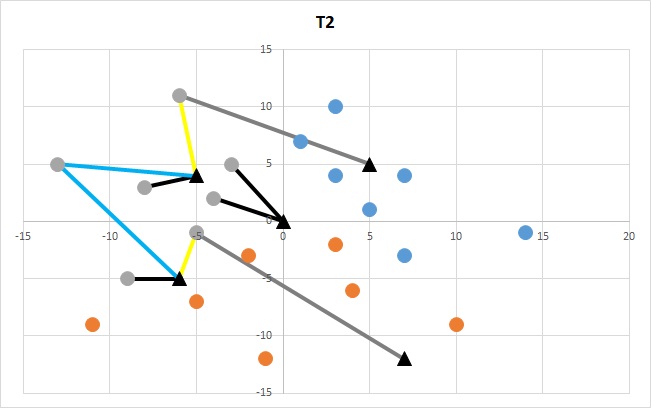}%
		\caption{Time Period 2}
	\end{subfigure}
	\\[0.3cm]
	\begin{subfigure}{0.5\textwidth}
		\centering
		\includegraphics[width=0.73\linewidth]{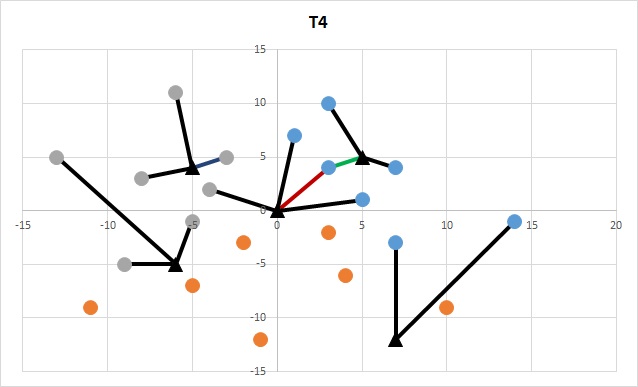}%
		\caption{Time Period 4}
	\end{subfigure}
	\begin{subfigure}{0.5\textwidth}
		\centering
		\includegraphics[width=0.73\linewidth]{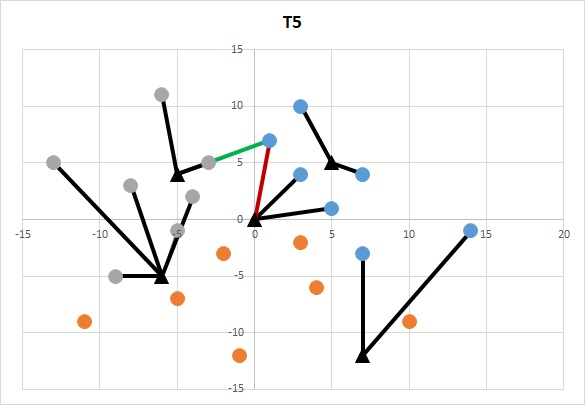}%
		\caption{Time Period 5}
	\end{subfigure}
	\\[0.3cm]
	\begin{subfigure}{0.5\textwidth}
		\centering
		\includegraphics[width=0.73\linewidth]{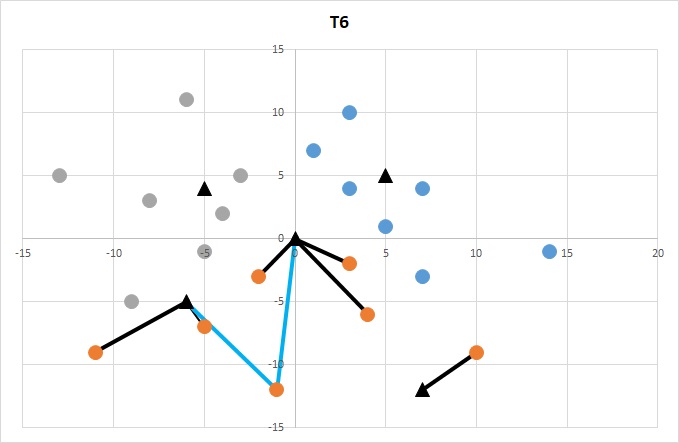}%
		\caption{Time Period 6}
	\end{subfigure}
	\begin{subfigure}{0.5\textwidth}
		\centering
		\includegraphics[width=0.73\linewidth]{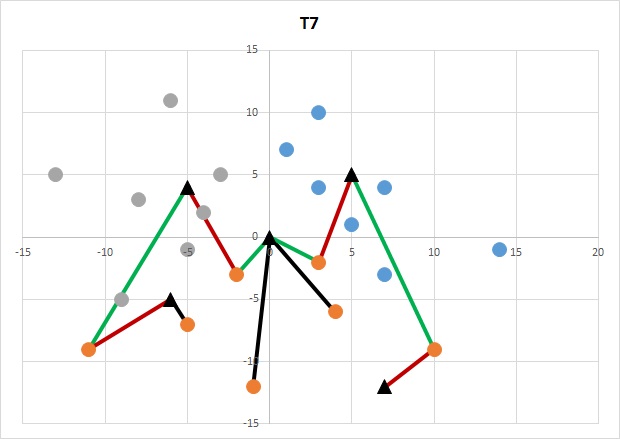}%
		\caption{Time Period 7}
	\end{subfigure}
	\caption{SEC instance: An optimal solution to the \mpmodel\ model.}
	\label{fig:M1_spicchi}
\end{sidewaysfigure}

\begin{figure}
	\centering
	\includegraphics[scale=0.75]{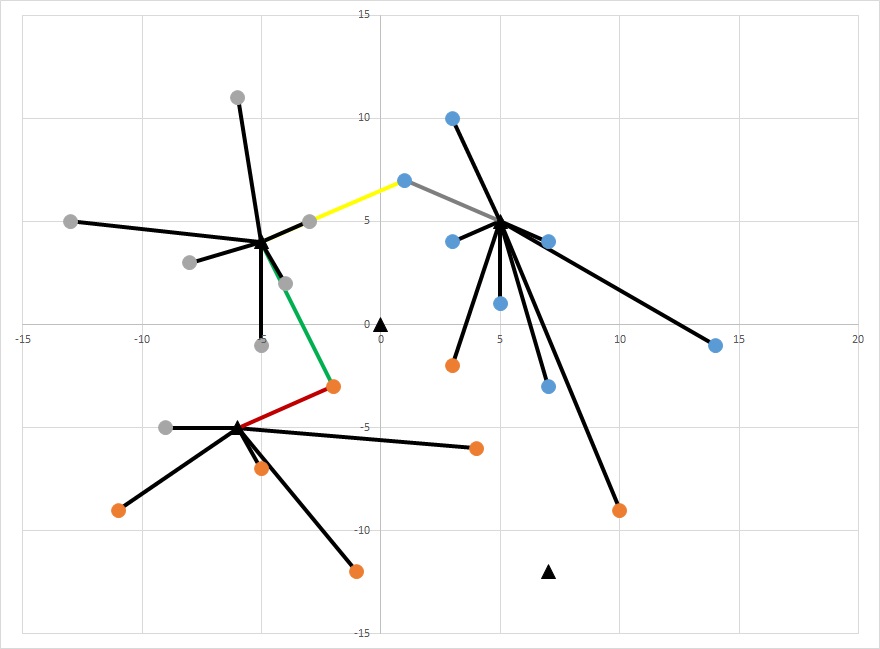}
		\caption{SEC instance: An optimal solution to the \spmodel\ model.}
		\label{fig:M4_spicchi}
\end{figure}

The optimal solutions found by the two MILP models for the SEC instance are shown in Figures~\ref{fig:M1_spicchi} and \ref{fig:M4_spicchi}.

Comparing the optimal solutions for the COR instance obtained by the \mpmodel\ and the \spmodel\ models (see Figures \ref{fig:M1_cerchi} and 
\ref{fig:M4_cerchi}, respectively), one can notice that in the former all the potential stations are open and the charging 
demand is assigned, in the majority of the cases, to the nearest station. The solution found by the \spmodel\ model 
 opens only three stations. This small example  highlights the limits of the latter model: it neglects that the 
charging demand is concentrated in few peak time periods, and, consequently,  underestimates the charging need in those time periods. In fact, most of the demand 
is assigned to the charging station located in the central position (coordinates (0,0)), but the chargers deployed there are not sufficient to serve all 
the EVs during the peak hours. Due to the lower number of stations opened (3 against 5) and the number of chargers approximately 40\% lower (30 against 
52), we observe that  17.04\% of customers  cannot be served by the solution to \spmodel.

Similar conclusions can be drawn observing the optimal solutions for the SEC instance produced by the \mpmodel\ and the \spmodel\ models (see 
Figures~\ref{fig:M1_spicchi} and \ref{fig:M4_spicchi}, respectively). As expected, there is no remarkable difference between the computational time on the COR and SEC instances. 
To keep a reasonable length of the paper, we  conducted the extensive experiments on the COR instances only.

\subsection{Computational results}\label{sec:CompRes}

\begin{figure}[H]
		\centering
	\includegraphics[scale=0.60]{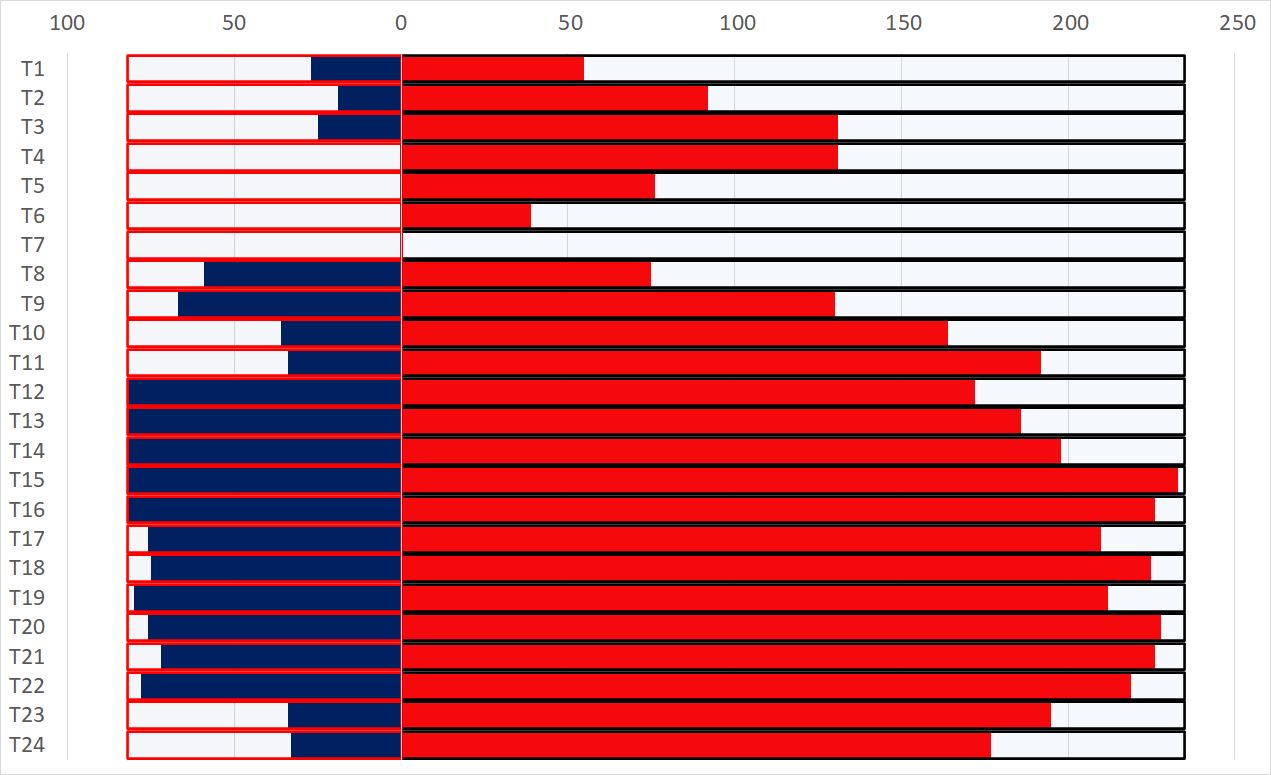}
	\caption{An optimal solution to the \mpmodel\ model for instance 200\textunderscore 30\textunderscore 30: Charging capacity installed (bordered bars) 
and demand assigned (solid bars).}
	\label{fig:M1_200_C_C}
\end{figure}

This section is devoted to the illustration and comment of the computational results. Before entering into the details of the results, we 
illustrate the solutions produced by the two MILP models for instance 200\textunderscore 30\textunderscore 30 \textbf and $\lambda = 0.50$. 
 Figure~\ref{fig:M1_200_C_C} depicts, for each time period, the charging capacity installed and the demand satisfied by
 an optimal solution to model \mpmodel\ for instance 
200\textunderscore 30\textunderscore 30. For each time period (vertical axis), the tornado diagram shows as bordered bars the number of chargers of each 
type deployed: the red bordered bars (left) are the fast chargers, whereas the black bordered bars (right) are the quick chargers. The solid bars represent 
the assignment of the charging demand. For each time period and type of charger, the bar indicates the total number of chargers assigned to 
EVs. Recall that quick chargers need multiple time periods to fully charge an EV. Hence, an EV assigned to a quick charger will use it for 
multiple consecutive time periods. From Figure~\ref{fig:M1_200_C_C}, one can notice that the demand assigned to each type of chargers in each time 
period does not violate the charging capacity deployed. Note also that in several time periods the demand approaches the capacity installed, and these two quantities are sometimes even equal (see the fast chargers in time periods 12 through 16).

\begin{figure}[H]
	\centering		
		\includegraphics[scale=0.60]{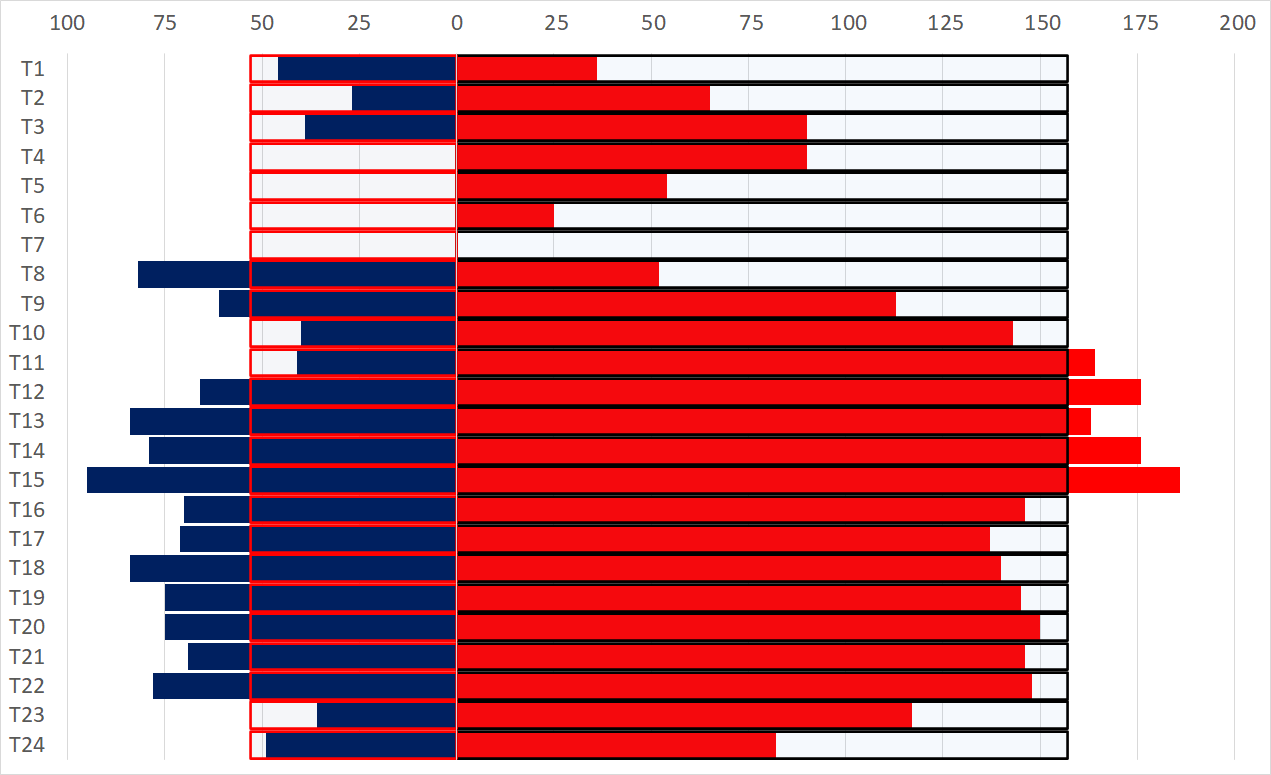}
		\caption{An optimal solution to the \spmodel\ model for instance 200\textunderscore 30\textunderscore 30: Charging capacity installed (bordered bars) and demand assigned (solid bars).}
		\label{fig:M4_200_C_C}		
	\end{figure}

The limits of the solution produced by the \spmodel\ model are evident from Figure~\ref{fig:M4_200_C_C}. 
The solution found by the \spmodel\ model  assigns, in several 
time periods, more EVs than the chargers actually available. For the fast chargers, this happens in time periods 8, 9, and from 12 to 22. On the 
other hand, for the quick chargers, this occurs in time periods 11 through 15. We observed this outcome 
for the majority of the instances tested. 
	
We now analyze more thoroughly the solutions produced by the two MILP models. To gain some insights about the two  
components of the objective functions, each instance is solved by each model for several values of the trade-off parameter $\lambda$. We 
tested the following values: $\lambda = 0.0001$ (maximum weight on the minimization of the total opening and installing costs), 0.25, 0.5, 0.75, and 0.9999 (maximum weight on the minimization of the average distance traveled by the EVs).

Table~\ref{t:solutions_1} provides, in the first three groups of columns, a summary of the charging capacity deployed by the solutions found 
by the two MILP models. For each group of instances and each model, 
Table~\ref{t:solutions_1} shows the average number of stations open (columns with header \textit{``Stations''}), as well as the average number of quick 
and fast chargers installed. For each value of $\lambda$, we reported in bold the average value of each of the former statistics for the \mpmodel\ model, 
along with the average deviation from the latter value for the \spmodel\ model.
The last group of three columns provides some statistics about the solution of the \spmodel\ model. In fact, the statistics refer to a modified solution obtained as follows. The deployed capacity remains unchanged. However, as the solution assigns the demand to open stations that may be  overloaded in some peak periods of time, we modified the assignment of the demand to the stations
with the goal of increasing the percentage of demand satisfied by the charging capacity deployed by the solution to the \spmodel\ model.
 
The procedure to modify the solution to the \spmodel\ model iteratively considers one time period at a time, from 1 to $T$, and, for a given time period $t$, examines each demand node $i$, from 1 to $I$. The procedure checks whether the demand $d_i^t$ of node $i$ could be served in time period $t$ according to the assignment indicated by the values of variables $x_{ijk}$. In this context, being served means that there is a number of vacant chargers $k$ in station $j$ greater than or equal to $d_i^t \cdot x_{ijk}$.

\begin{table}
	\centering
	\resizebox*{0.90\textwidth}{!}{
		\begin{tabular}[h!]{cr|rr|rr|rr||rrr}
			\toprule
			& &  \multicolumn{2}{|c|}{\textbf{Stations}}&\multicolumn{2}{|c|}{\textbf{Quick}}&\multicolumn{2}{|c||}{\textbf{Fast}} 
&\multicolumn{3}{|c}{\textbf{\spmodel}}\\
			$\bm{\lambda}$ & $\bm{I$}  & \textbf{\mpmodel} & \textbf{\spmodel} & \textbf{\mpmodel} & \textbf{\spmodel} & \textbf{\mpmodel} & 
\textbf{\spmodel} & \textbf{Reall\%} & \textbf{Lost\%} & \textbf{Max Lost\%}\\
			\midrule
			\multirow{6}{*}{\textbf{0.0001}}
&	50	&	4.33	&	2.67	&	53.33	&	33.00	&	23.33	&	15.00 &	6.37\%	&	21.00\%	&	60.23\%	\\
&	100	&	8.20	&	5.00	&	101.07	&	60.53	&	45.87	&	30.40	&	9.82\%	&	21.07\%	&	58.33\%	\\
&	150	&	10.87	&	7.53	&	141.40	&	91.93	&	60.73	&	45.93	&	10.95\%	&	20.88\%	&	58.03\%	\\
&	200	&	14.64	&	9.71	&	186.00	&	126.50	&	94.14	&	61.00	&	10.41\%	&	21.10\%	&	58.60\%	\\
&	250	&	17.75	&	12.00	&	239.42	&	153.93	&	115.08	&	77.36	&	11.24\%	&	21.28\%	&	58.19\%	\\
&	500	&	31.25	&	21.64	&	418.56	&	300.18	&	198.78	&	154.82	&	13.24\%	&	21.82\%	&	55.89\%	\\ 		
			\midrule
			\textbf{Average} & &  \textbf{12.90} & \textbf{-28.93\%} & \textbf{171.01} & \textbf{-30.32\%} & \textbf{80.46} & \textbf{-25.87\%} &	
\textbf{10.19\%}	&	\textbf{21.16\%}	&	\textbf{58.32\%}\\
		\midrule
			\multirow{6}{*}{\textbf{0.25}}
&	50	&		4.80	&	3.80	&	55.47	&	35.60	&	22.80	&	14.13	&	12.23\%	&	21.15\%	&	55.90\%	\\
&	100	&		8.13	&	5.33	&	95.00	&	61.73	&	47.07	&	29.93	&	12.12\%	&	20.90\%	&	58.43\%	\\
&	150	&		11.07	&	7.73	&	137.00	&	90.13	&	68.57	&	46.07	&	11.66\%	&	20.98\%	&	57.90\%	\\
&	200	&		14.14	&	9.79	&	170.14	&	120.64	&	97.93	&	62.14	&	12.17\%	&	21.11\%	&	57.73\%	\\
&	250	&		17.50	&	11.79	&	227.58	&	143.21	&	117.75	&	79.79	&	12.54\%	&	21.35\%	&	54.93\%	\\
&	500	&		30.38	&	21.27	&	393.44	&	287.27	&	229.75	&	159.00	&	14.56\%	&	21.79\%	&	54.95\%	\\
			\midrule
\textbf{Average} & &  \textbf{12.82} & \textbf{-26.74\%} & \textbf{162.39} & \textbf{-29.14\%} & \textbf{85.00} & \textbf{-28.74\%} &	\textbf{12.46\%}	
&	\textbf{21.19\%}	&	\textbf{56.73\%}		\\
		\midrule
			\multirow{6}{*}{\textbf{0.50}}
&	50	&	6.93	&	6.60	&	62.67	&	36.60	&	23.20	&	14.40	&	18.38\%	&	19.64\%	&	54.65\%	\\
&	100	&	8.67	&	7.60	&	97.73	&	73.67	&	46.20	&	27.00	&	14.02\%	&	20.56\%	&	62.05\%	\\
&	150	&	11.07	&	8.80	&	131.86	&	97.87	&	69.93	&	44.53	&	13.55\%	&	20.58\%	&	59.55\%	\\
&	200	&	14.21	&	10.71	&	166.71	&	126.29	&	98.86	&	60.86	&	13.93\%	&	20.82\%	&	58.03\%	\\
&	250	&	17.75	&	12.50	&	219.67	&	150.29	&	119.58	&	78.00	&	14.41\%	&	21.06\%	&	57.80\%	\\
&	500	&	30.88	&	21.55	&	441.50	&	285.73	&	230.50	&	158.82 &	16.00\%	&	21.69\%	&	54.80\%	\\
\midrule
\textbf{Average} & &  \textbf{13.44} & \textbf{-19.64\%} & \textbf{163.51} & \textbf{-26.20\%} & \textbf{85.68} & \textbf{-30.81\%} &	\textbf{14.86\%}	
&	\textbf{20.68\%}	&	\textbf{57.95\%}	\\
		\midrule		
\multirow{6}{*}{\textbf{0.75}}
&	50	&	12.13	&	12.20	&	78.60	&	37.47	&	31.07	&	15.20	&	26.82\%	&	17.09\%	&	51.54\%	\\
&	100	&	12.20	&	11.40	&	120.60	&	70.27	&	50.20	&	28.53	&	19.83\%	&	19.71\%	&	60.30\%	\\
&	150	&	13.64	&	12.80	&	152.93	&	103.73	&	73.14	&	43.00	&	18.00\%	&	20.06\%	&	60.58\%	\\
&	200	&	15.21	&	13.50	&	173.36	&	137.07	&	97.14	&	58.14	&	16.34\%	&	20.59\%	&	59.90\%	\\
&	250	&	18.00	&	14.79	&	223.25	&	156.71	&	120.58	&	76.43	&	16.13\%	&	20.88\%	&	58.69\%	\\
&	500	&	30.63	&	21.91	&	428.25	&	285.36	&	233.75	&	157.91	&	17.78\%	&	21.66\%	&	54.67\%	\\
			\midrule
\textbf{Average} & &  \textbf{15.77} & \textbf{-10.69\%} & \textbf{175.14} & \textbf{-29.15\%} & \textbf{88.72} & \textbf{-33.95\%} &	\textbf{19.02\%}	
&	\textbf{19.90\%}	&	\textbf{57.71\%}	\\
		\midrule
\multirow{6}{*}{\textbf{0.9999}}
&	50	&	21.20	&	21.20	&	96.87	&	40.53	&	41.27	&	17.00	&	31.11\%	&	12.33\%	&	47.73\%	\\
&	100	&	27.20	&	27.20	&	160.20	&	74.67	&	77.60	&	31.40	&	26.57\%	&	16.04\%	&	52.02\%	\\
&	150	&	29.71	&	28.40	&	212.93	&	104.33	&	113.86	&	46.87	&	24.51\%	&	16.62\%	&	57.10\%	\\
&	200	&	30.93	&	30.07	&	245.50	&	141.43	&	146.00	&	61.86	&	23.28\%	&	17.41\%	&	57.54\%	\\
&	250	&	33.58	&	30.79	&	276.00	&	165.29	&	158.29	&	79.43	&	22.05\%	&	17.92\%	&	57.05\%	\\
&	500	&	38.75	&	34.55	&	429.38	&	285.45	&	375.25	&	165.27	 &	19.59\%	&	19.90\%	&	52.74\%	\\
			\midrule
\textbf{Average} & &  \textbf{29.33} & \textbf{-3.25\%} & \textbf{218.95} & \textbf{-41.67\%} & \textbf{132.99} & \textbf{-53.23\%} &	\textbf{24.80\%}	
&	\textbf{16.53\%}	&	\textbf{54.02\%}	\\
\bottomrule		
		\end{tabular}
	}
		\caption{\mpmodel\ vs. \spmodel\ models: An analysis of the solutions found.}
		\label{t:solutions_1}
\end{table}

If such demand cannot be completely served, the unserved demand is reallocated among the vacant chargers of any type different 
from $k$ available at the same station $j$, if any. 
Then, the procedure attempts to reallocate the remaining unserved demand among the other stations. If there are vacant chargers among multiple stations, 
priority is given to the one nearest to $j$.  If at the  station there are vacant 
chargers of multiple types, priority is given to the type that has the largest number of vacant units. 
Let $0 \leq \bar{d}_i^t \leq d_i^t$ be the demand arising in node $i$ at time period $t$ that the procedure has reallocated. 
Eventually, the procedure computes the fraction of the total demand that has been reallocated, in percentage, producing statistic \textit{``Reall\%''}. For 
each instance, the fraction of the total demand reallocated is computed as $100 \cdot \frac{\sum_{t \in \mathcal{T}}\sum_{i \in \mathcal{I}} 
\bar{d}_i^t}{\sum_{t \in \mathcal{T}}\sum_{i \in \mathcal{I}} d_i^t}$.

If, after the reallocation procedure, a portion of the demand is not  served, the fraction of the total demand that remains unserved is computed, in percentage, 
producing statistic \textit{``Lost\%''}. For each instance, the latter is computed as $100 \cdot \frac{\sum_{t \in \mathcal{T}}\sum_{i \in \mathcal{I}} 
\check{d}_i^t}{\sum_{t \in \mathcal{T}}\sum_{i \in \mathcal{I}} d_i^t}$, where $0 \leq \check{d}_i^t \leq d_i^t$ is the demand arising in node $i$ at time 
period $t$ that is not served. Finally, the procedure computes statistic \textit{``Max Lost\%''} as the maximum fraction of demand not served across all 
time periods. The latter is computed for each instance as $100 \cdot \max\limits_{t \in \mathcal{T}} \left\{ \frac{\sum_{i \in \mathcal{I}} 
\check{d}_i^t}{\sum_{i \in \mathcal{I}} d_i^t} \right\}$.

From Table~\ref{t:solutions_1} we can gain the following insights:

\begin{itemize}
	\item[\checkmark] the charging capacity deployed in the solutions to the \spmodel\ model is significantly smaller than the capacity installed 
according to the \mpmodel\ model, both in terms of stations open and chargers installed (see statistics ``Stations'', ``Quick'', and ``Fast'');
	\item[\checkmark] reducing the weight of the  opening and installing costs (i.e., increasing the value of $\lambda$), the average deviation 
between the solutions to the two models in terms of stations open decreases steadily, from -28.93\% for $\lambda = 0.0001$ to -3.25\% for $\lambda = 
0.9999$. Nevertheless, the average number of chargers installed (of each type) in the solutions found by the \spmodel\ model is always remarkably smaller 
compared to those produced by the \mpmodel\ model;
	\item[\checkmark] given a value of $\lambda$, for both models the greater the number of demand nodes, the greater the number of stations open and 
chargers installed (see statistics ``Stations'', ``Quick'', and ``Fast'');
	\item[\checkmark] the larger the value of $\lambda$, the larger the values of ``Stations'', ``Quick'', and ``Fast''. This is an expected outcome, as 
more importance is given to the average distance term in the objective functions. In other words, the reduction of the average distance traveled can only 
be achieved by increasing the number of stations opened and chargers deployed;
	\item[\checkmark] due to the  cheaper installation cost, both models install a larger number of quick compared to  fast chargers. 
\end{itemize}

Before entering into the details of the statistics computed to measure the limits of the \spmodel\ model, it is worth pointing out that 
in every solution found by the latter model, a part of the demand was reallocated, and a part was lost. The main insights we 
can gain from the three rightmost columns of Table~\ref{t:solutions_1} are the following:

\begin{itemize}
	\item[\checkmark] ``Reall\%'' takes, on average, large values. It ranges from 6.37\% (see $\lambda = 0.0001$ and $I = 50$) to 31.11\% (see $\lambda = 
0.9999$ and $I = 50$);
	\item[\checkmark] ``Lost\%'' takes, on average, large values as well. It ranges from 12.33\% (see $\lambda = 0.9999$ and $I = 50$) to 21.82\% (see 
$\lambda = 0.0001$ and $I = 500$);
	\item[\checkmark] ``Max Lost\%'' takes, on average, extremely large values, always larger than 47\%. 
\end{itemize}

\begin{figure}
	\begin{subfigure}{0.5\textwidth}
		\centering
		\includegraphics[width=0.9\linewidth]{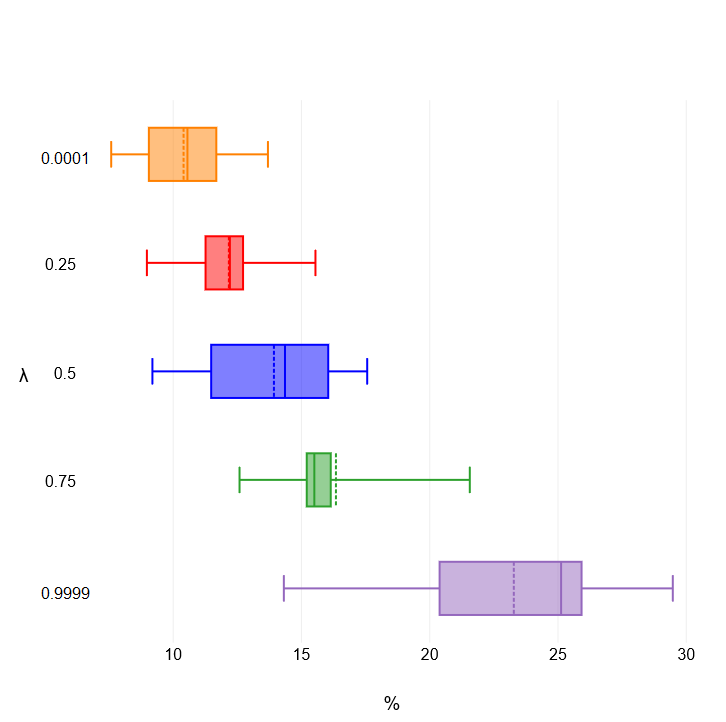}%
		\caption{Demand reallocated (``Reall\%'')}
	\end{subfigure}
	\begin{subfigure}{0.5\textwidth}
		\centering
		\includegraphics[width=0.9\linewidth]{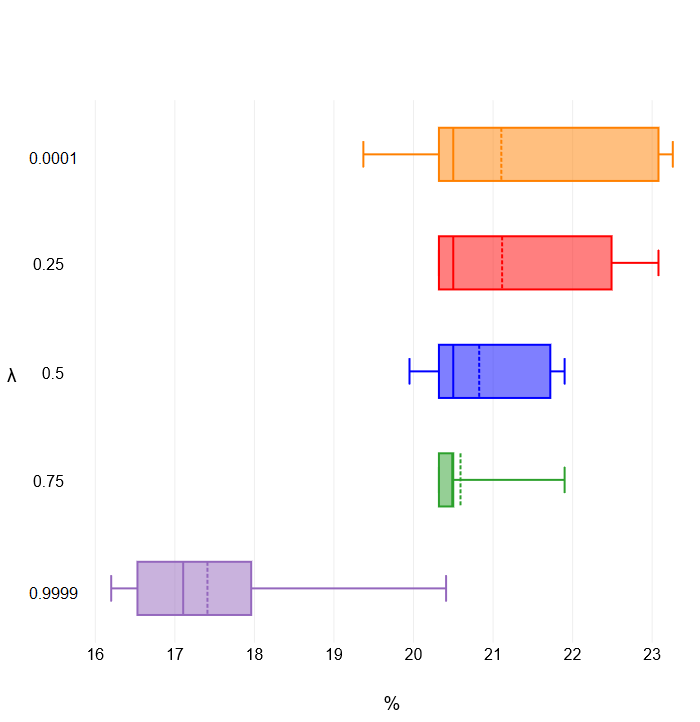}%
		\caption{Unserved demand (``Lost\%'')}
	\end{subfigure}
	\\[0.3cm]
		\begin{subfigure}{1.\textwidth}
		\centering
		\includegraphics[width=0.5\linewidth]{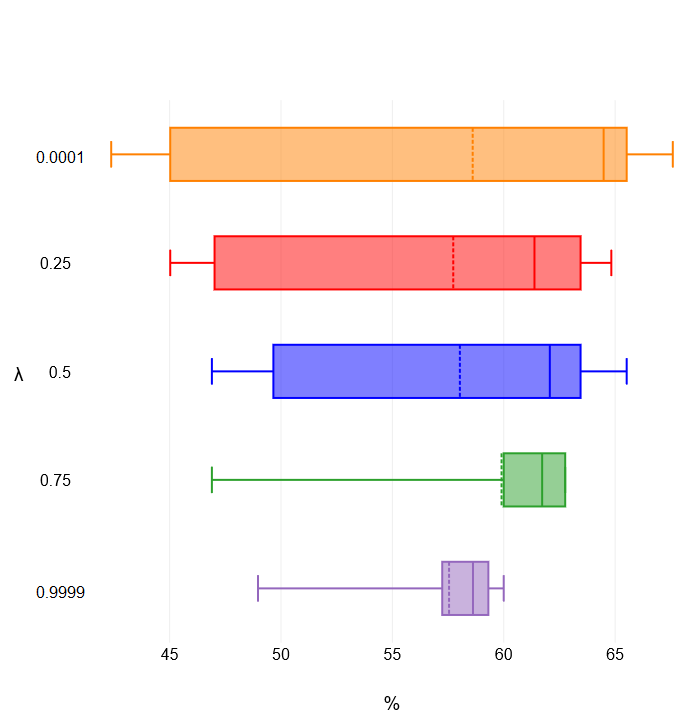}%
		\caption{Maximum fraction of unserved demand   across all time periods (``Max Lost\%'')}
	\end{subfigure}
		\caption{\spmodel\ model: Box-and-wisker plots showing the distribution of the demand reallocated and lost ($I=200$).}		
		\label{fig:Repo}%
\end{figure}

The statistics confirm that the \spmodel\ model is not capable of  capturing  the characteristics of the problem and tends to 
underestimate the charging capacity to deploy.

Further insights that can be obtained from Table~\ref{t:solutions_1} on the limits of the \spmodel\ model are as follows:

\begin{itemize}
	\item[\checkmark] for values of $\lambda$ smaller than or equal to 0.25, the average value of ``Reall\%'' tends to increase with the number of demand 
nodes;
	\item[\checkmark] for values of $\lambda$ greater than or equal to 0.75, the average value of ``Reall\%'' tends to decrease with the number of demand 
nodes;
	\item[\checkmark] the larger the value of $\lambda$, the larger the average value of ``Reall\%'', and the smaller tends to be the value of 
``Lost\%''.  This behavior can be explained by observing that increasing the value of $\lambda$, the number of stations opened and chargers installed 
increases as well, making it easier to find vacant chargers, and thereby reducing the unserved demand.
of ``Max Lost\%'' slightly decreases as the value of $\lambda"$ increases.
\end{itemize}

The box-and-whisker plots depicted in Figure~\ref{fig:Repo} show the distribution of the values of the three statistics computed to determine 
the  reallocated demand, and the unserved demand, for all the instances with $I=200$  solved for different values of $\lambda$. The box-and-whisker 
plots confirm the insights  previously drawn. When the main term in the objective function of the \spmodel\ model is the  opening and installing cost 
- i.e, for small values of $\lambda$ - the charging capacity installed is small and little can be done to reallocate the unserved demand. 
As a consequence, large 
percentages of the charging demand are unserved. By increasing the weight given to the average distance traveled - i.e., for large values of $\lambda$ - the 
charging capacity installed increases. Consequently, the percentage of the demand that can be reallocated increases, and, thereby, the percentage of the 
demand that is unserved becomes smaller. Nevertheless, the latter percentages always remain quite large (the values are always greater than 16\%, 
and often greater than 20\%). The performance is even worse if we analyze the distribution of ``Max Lost\%''. Its average value ranges from approximately 56\% to roughly 60\% (see the dotted lines inside the boxes). 
 Similar conclusions can be drawn observing the results obtained 
when solving the instances with a different number of demand nodes. For the sake of readability, such results are not reported here.

\begin{figure}
	\begin{subfigure}{0.5\textwidth}
		\centering
		\includegraphics[width=0.9\linewidth]{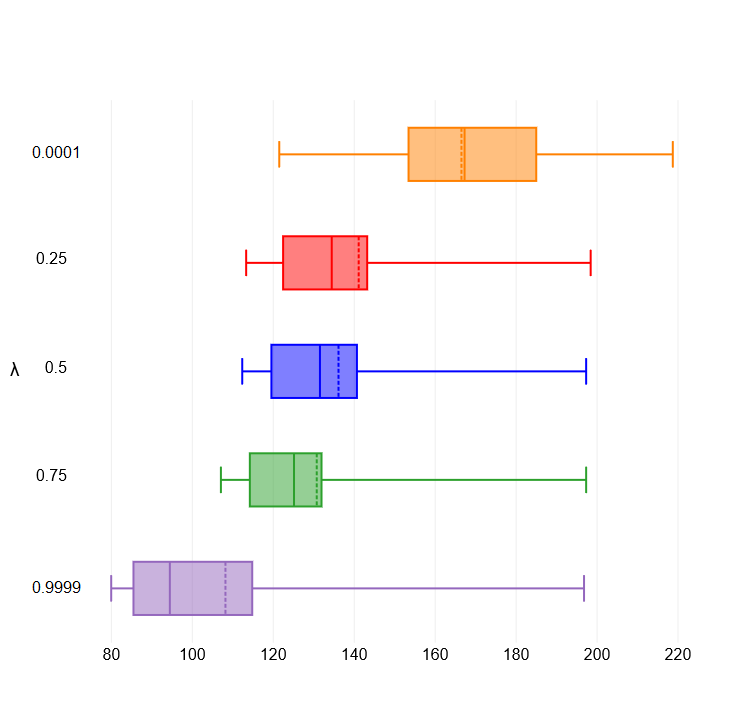}
		\caption{Distribution of average distance traveled.}
	\end{subfigure}
			\begin{subfigure}{0.5\textwidth}
		\centering
		\includegraphics[width=0.9\linewidth]{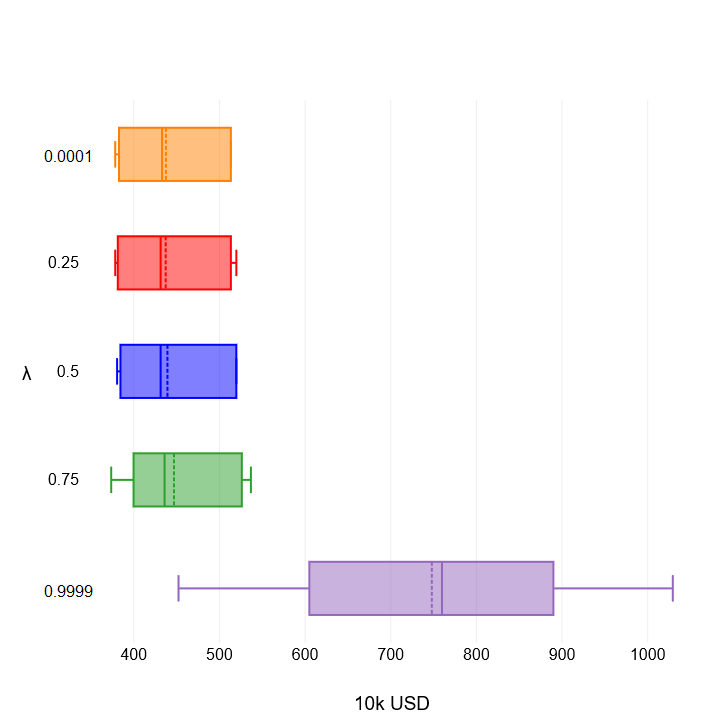}
		\caption{Distribution of total opening and installing costs.}
	\end{subfigure}
	\caption{\mpmodel\ model: Box-and-wisker plots showing the distribution of the average distance and the cost ($I=200$).}\label{fig:costs_distance}
\end{figure}

The results discussed above clearly show that the solutions found by the \spmodel\ model, if implemented, would lead to a 
very poor quality of service provided to the EV drivers. 
Moreover, the results on demand reallocation imply that the objective function of \spmodel\ would underestimate the total travel distance covered by the EV drivers to reach a free charging station.

The following analysis is focused only on the solutions 
of the \mpmodel\ model. 
Figure~\ref{fig:costs_distance} illustrates the distributions of the values of the average distance traveled (first term in the objective function) and the  opening and installing cost (second term), for all the instances with $I=200$  solved for different 
values of $\lambda$.

From Figure~\ref{fig:costs_distance}, we can draw the following main insights:

\begin{itemize}
\item[\checkmark] as expected, the larger the value of $\lambda$, the smaller the average distance traveled by an EV to reach the assigned charger;
\item[\checkmark] besides its largest value, increasing the value of $\lambda$ produces only slight increases in the values of the total cost.
\end{itemize}

In fact, we expected a sharper increase of the values of the total cost when less importance is given to the second term of the objective function. On the contrary, and neglecting the extreme case with $\lambda = 0.9999$, when the value of $\lambda$ is increased the solutions obtained by the \mpmodel\ model  significantly improved the average traveling distance, at the cost of only a small deterioration of the total opening and installing cost.

We conclude our analysis by considering the computational burden required to solve each MILP model. Recall that each instance is solved with a time limit 
of 3,600 seconds. 
Table~\ref{t:comp_2} summarizes the computational performance of the two MILP models. For each value of $\lambda$, the instances are clustered in groups 
according to the number of potential locations $J$. For each group of instances and each model, Table~\ref{t:comp_2} provides the average CPU time (in 
seconds) spent to find the optimal (or best) solution (columns with header \textit{``CPU Time (secs.)''}), the average optimality gap (\textit{``Gap\%''}) 
and the worst optimality gap (\textit{``Max Gap\%''}).

The main insights that we can gain from Table~\ref{t:comp_2} are as follows:

\begin{itemize}
	\item[\checkmark] the solution to the \mpmodel\ model is, in general, more computationally expensive compared to the \spmodel\ model (see the average values of ``CPU Time (secs.)'', and the average values of ``Gap\%'');
	\item[\checkmark] the optimality gaps, for both models, are on average very small. In the majority of the cases, the solver found an optimal 
solution, or a solution very close to the optimum;
	\item[\checkmark] the worst gaps, for both models, are also very small. In only few instances, statistic ``Max Gap\%'' took a value greater than 1\%;
	\item[\checkmark] as expected, for a given value of $\lambda$, computing times for both models increase with the number of potential locations;	
	\item[\checkmark] the computational burden required to solve each MILP model decreases as the value of $\lambda$ increases.
\end{itemize}

In summary, while  the \mpmodel\ model  requires on average more computational time than the \spmodel\ model, the additional time 
needed by the \mpmodel\ model is marginally small.

\begin{table}
	\centering
	\resizebox*{0.65\textwidth}{!}{
		\begin{tabular}[h!]{cr|rr|rr|rr}
			\toprule
			& &  \multicolumn{2}{|c|}{\textbf{CPU Time (secs.)}}&\multicolumn{2}{|c|}{\textbf{Gap\%}}&\multicolumn{2}{|c}{\textbf{Max Gap\%}} \\
			$\bm{\lambda}$ & $\bm{J}$  & \textbf{\mpmodel} & \textbf{\spmodel} & \textbf{\mpmodel} & \textbf{\spmodel} & \textbf{\mpmodel} & 
\textbf{\spmodel} \\
			\midrule
			
			\multirow{5}{*}{\textbf{0.0001}}
&	10	&	1,456.48	&	820.60	&	0.13\%	&	0.09\%	&	0.45\%	&	0.76\%	\\
&	20	&	2,775.72	&	2,404.48	&	0.61\%	&	0.46\%	&	1.85\%	&	1.87\%	\\
&	30	&	3,308.19	&	2,759.82	&	1.03\%	&	1.08\%	&	2.90\%	&	3.69\%	\\
&	40	&	3,515.42	&	3,142.17	&	0.94\%	&	1.23\%	&	2.90\%	&	3.69\%	\\
&	50	&	3,513.11	&	3,330.17	&	1.08\%	&	1.43\%	&	2.90\%	&	3.69\%	\\
\midrule
\textbf{Average} & &  \textbf{3,037.11}	&	\textbf{2,568.85}	&	\textbf{0.81\%}	&	\textbf{0.90\%}	&	& \\ 
\midrule			
			\multirow{5}{*}{\textbf{0.25}}
&	10	&	660.83	&	259.72	&	0.03\%	&	0.01\%	&	0.17\%	&	0.11\%	\\
&	20	&	2,133.02	&	922.11	&	0.12\%	&	0.04\%	&	0.60\%	&	0.37\%	\\
&	30	&	2,923.03	&	1,320.37	&	0.31\%	&	0.13\%	&	0.92\%	&	0.82\%	\\
&	40	&	3,156.36	&	1,790.20	&	0.45\%	&	0.19\%	&	1.19\%	&	1.07\%	\\
&	50	&	3,406.24	&	2,120.10	&	0.55\%	&	0.18\%	&	1.77\%	&	0.84\%	\\
\midrule
\textbf{Average} & &  \textbf{2,614.44}	&	\textbf{1,335.04}	&	\textbf{0.32\%}	&	\textbf{0.11\%}	&	& \\ 
\midrule		
\multirow{5}{*}{\textbf{0.50}}
&	10	&	263.90	&	2.43	&	0.00\%	&	0.00\%	&	0.00\%	&	0.00\%	\\
&	20	&	1,135.25	&	338.01	&	0.01\%	&	0.01\%	&	0.06\%	&	0.17\%	\\
&	30	&	2,544.90	&	1,334.12	&	0.19\%	&	0.07\%	&	0.67\%	&	0.66\%	\\
&	40	&	2,739.82	&	2,051.00	&	0.34\%	&	0.09\%	&	1.06\%	&	0.58\%	\\
&	50	&	2,978.22	&	2,145.28	&	0.97\%	&	0.10\%	&	6.85\%	&	0.49\%	\\
\midrule
\textbf{Average} & &  \textbf{2,094.62}	&	\textbf{1,238.01}	&	\textbf{0.34\%}	&	\textbf{0.06\%}	&	& \\ 
\midrule					
\multirow{5}{*}{\textbf{0.75}}
&	10	&	6.45	&	2.51	&	0.00\%	&	0.00\%	&	0.00\%	&	0.00\%	\\
&	20	&	376.49	&	1,183.04	&	0.00\%	&	0.01\%	&	0.02\%	&	0.13\%	\\
&	30	&	1,678.43	&	1,984.19	&	0.04\%	&	0.05\%	&	0.41\%	&	0.22\%	\\
&	40	&	2,449.72	&	2,298.07	&	0.06\%	&	0.07\%	&	0.37\%	&	0.24\%	\\
&	50	&	2,926.84	&	2,312.06	&	0.40\%	&	0.06\%	&	4.28\%	&	0.22\%	\\
\midrule
\textbf{Average} & &  \textbf{1,648.46}	&	\textbf{1,629.29}	&	\textbf{0.11\%}	&	\textbf{0.04\%}	&	& \\ 
\midrule					
\multirow{5}{*}{\textbf{0.9999}}
&	10	&	2.80	&	1.97	&	0.00\%	&	0.00\%	&	0.00\%	&	0.00\%	\\
&	20	&	966.18	&	1,488.73	&	0.00\%	&	0.00\%	&	0.00\%	&	0.00\%	\\
&	30	&	1,575.65	&	2,581.28	&	0.00\%	&	0.00\%	&	0.00\%	&	0.00\%	\\
&	40	&	1,864.13	&	3,243.86	&	0.00\%	&	0.00\%	&	0.00\%	&	0.00\%	\\
&	50	&	2,569.15	&	2,957.01	&	0.00\%	&	0.00\%	&	0.00\%	&	0.00\%	\\
\midrule
\textbf{Average} & &  \textbf{1,528.78}	&	\textbf{2,152.78}	&	\textbf{0.00\%}	&	\textbf{0.00\%}	&	& \\ 
\bottomrule		
		\end{tabular}
	}
		\caption{\mpmodel\ vs. \spmodel\ models: A summary of CPU times and optimality gaps.} \label{t:comp_2}
\end{table}

\section {Conclusions} \label{sec:conc}
In this paper, we studied the role of temporal and spatial  distributions of charging demand in determining an optimal location of charging stations for electrical vehicles in an urban setting. 
This is an application context where the daily demand is well-known to be very dynamic and concentrated at some peak hours, and where the demand pattern is known to depend also upon the city zone. 

To highlight the need of considering explicitly  
the daily demand patterns,
 we presented a multi-period optimization model, that captures the variability over time of the demand,
  and compared it with a single-period optimization model.
  By means of a worst-case analysis, we  theoretically proved that the single-period model may  produce solutions where a large portion of the demand 
cannot be served. 
Extensive computational experiments confirm the limits of the single-period model.  The goal of the optimization models is to  balance
 two objectives: the total cost of deploying the infrastructure and the average distance traveled 
by the customers to 
reach a charging station.

The limits pointed out for the single-period model go beyond the specific application considered in this paper, and suggest the importance of incorporating 
time-dependency in location decisions when the demand fluctuations are remarkable during the planning horizon.

Future developments of the  research may concern the objective function. Since the two objectives considered are not homogeneous, a thorough 
analysis of the trade-off between infrastructure cost and drivers traveling distance may be of interest to a decision-maker. Moreover, we expect that 
replacing the average traveling distance with some equity measures may produce solutions that are 
more satisfactory for the customers, at the cost of a little increase in the infrastructure cost.
Finally, to solve larger instances, in particular when an equity measure is considered, a heuristic approach would deserve to be studied.

{\subsection*{Acknowledgements} This study has been made in the framework of the MoSoRe@UniBS 
(Infrastrutture e servizi per la Mobilità Sostenibile e Resiliente) Project 2020-2022 of Lombardy Region, 
Italy (Call-Hub ID 1180965; bit.ly/2Xh2Nfr, https://ricerca2.unibs.it/?page id=8548).

\bibliographystyle{apalike}
\bibliography{literature}

\end{document}